\documentclass[twocolumn]{autart}    % Enable this line and disable the 
\usepackage{graphicx}
\usepackage{amsmath}
\usepackage{amssymb}
\usepackage{times,mathptm}
\usepackage{amsfonts}
\usepackage{enumerate}
\usepackage{longtable}
\usepackage{graphicx}
\usepackage{theorem}
\usepackage{multirow}
\usepackage{mdframed}
\usepackage{float}    
\usepackage{array} 
\usepackage[square,sort,comma,numbers]{natbib}

\newtheorem{proof}{Proof} [section]
\newtheorem{remark}{Remark} [section]
\newtheorem{lemma}{Lemma} [section]

\newtheorem{assumption}{Assumption} [section]
\newtheorem{theorem}{Theorem} [section]

\begin{document}

\begin{frontmatter}

\title{Useful redundancy in parameter and time delay estimation for continuous-time models}

\thanks[footnoteinfo]{The material in this paper was not presented at any conference.}

\author[First]{Huong Ha}\ead{huong.ha@uon.edu.au},
\author[First]{James S. Welsh}\ead{james.welsh@newcastle.edu.au},
\author[Second]{Mazen Alamir}\ead{mazen.alamir@grenoble-inp.fr},

\address[First]{School of Electrical Engineering and Computer Science, The University of Newcastle, Australia}
\address[Second]{Gipsa-lab, Control Systems Department. University of Grenoble, France}


\begin{abstract}                          % Abstract of not more than 200 words.
This paper proposes an algorithm to estimate the parameters, including time delay, of continuous time systems based on instrumental variable identification methods. To overcome the multiple local minima of the cost function associated with the estimation of a time delay system, we utilise the useful redundancy technique. Specifically, the cost function is filtered through a set of low-pass filters to improve convexity with the useful redundancy technique exploited to achieve convergence to the global minimum of the optimization problem. Numerical examples are presented to demonstrate the effectiveness of the proposed algorithm.
\end{abstract}
\end{frontmatter}

\vspace{-0.8cm}
\section{Introduction}
\vspace{-0.5cm}
The goal of system identification is to estimate the parameters of a model in order to analyse, simulate and/or control a system. In the time domain, there are two typical approaches to identify a system, i.e. discrete-time (DT) identification and continuous-time (CT) identification. For several decades, DT identification has been dominant due to the strong development of the digital computer. More recently, estimation using continuous-time identification methods has received much attention due to advantages such as providing insights to the physical system and being independent of the sampling time \cite{Garnier2015} \cite{Garnier2013} \cite{Rao1990}. For example, with irregular sampling time, the DT model becomes time-varying; hence the DT system identification problem becomes more difficult while for CT identification, the system is still time-invariant. In reality, irregular sampling occurs in many cases, e.g. when the sampling is event-triggered, when the measurement is manual and also in the case of missing data \cite{Astrom03}.

In some CT identification problems, one wants to estimate the system parameters and any unknown time delay, as there exist many practical examples including, chemical processes, economic systems, biological systems, that possess time delays. It is important to estimate the delay accurately since a poor estimate can result in poor model order selection and inaccurate estimates of the system parameters. 

There are many approaches to estimate a system time delay~\cite{Bjorklund2003}. A simple approach is to consider the impulse response data, e.g. estimate the time delay by finding where the impulse response becomes nonzero \cite{Carlemalm1999} or by noting the delay where the correlation between input and output is maximum \cite{Carlemalm1999} \cite{Carter87}. Another approach is to model the delay by a rational polynomial transfer function using a Pad\'{e} or similar approximation and then estimate the time delay as part of the system parameters \cite{Agarwal87} \cite{Ahmed2006} \cite{Gawthrop1985}. In \cite{Yang2007} \cite{Baysse2011} \cite{Baysse2012}, the time delay and system parameters of a Multiple Input Single Output CT system are estimated in a separable way using an iterative global nonlinear least-squares or instrumental variable method. 
%In \cite{Ni2010}, the cross correlation between the continuous wavelet transform coefficient of the input and output data is computed in order to estimate the unknown time delays for MIMO systems. 

Recently, a method of estimating the parameters and time delay of CT systems has been suggested in \cite{Chen2015} which is based on a gradient technique. The parameters and the time delay are estimated separately, i.e. when one is estimated, the other is fixed which is then repeated in an iterative manner. In this approach, the Simplified Refined Instrumental Variable (SRIVC) method is used to estimate the parameters whilst the time delay is estimated using the Gauss-Newton method. In addition, due to the effects of multiple minima in the cost function to be minimized for the time delay \cite{Kaminskas1979} \cite{Pupeikis1989} \cite{Bjorklund1991}, a low-pass filter is employed to increase convexity. As shown in \cite{Bjorklund1991} \cite{Ferretti:96} \cite{Diego2016}, a suitable low-pass filtering operation on the estimation data can help to extend the global convergence region of the cost function, hence improve the accuracy of the time delay estimate.

In this paper, we adopt the idea of using a low-pass filter where instead of using only one filter, we suggest to use multiple low-pass filters and incorporate the useful redundancy technique \cite{Alamir2008} \cite{Alamir2009}. Useful redundancy is a technique that was developed to avoid local minima when solving a nonlinear inverse problem. The concept is to generate a family of cost functions that have different local minima, but share the same global minimum with the original cost function of the optimization problem. Whenever the algorithm is trapped in a minimum, the solver path is switched to another solver path in a way such that the minimum found using the new solver path corresponds to a decrease in the original cost function. This allows the algorithm to cross local minima and converge to the global minimum, hence improving the accuracy of the estimated parameters. For the algorithm described in this paper, the multiple cost functions are generated by filtering the original time delay cost function through a number of low-pass filters with different cut-off frequencies that span the system bandwidth.

The paper is organized as follows. Section 2 describes the model setting and Section 3 recalls both the SRIVC algorithm and the SRIVC-based time delay estimation using a low-pass filter. Section 4 formulates the algorithm of the new method and provides analysis on the effectiveness of the proposed method. Section 5 presents numerical experiments and results for both regular and irregular sampling schemes. Finally, the conclusion will be drawn in Section 6.

%%%%%%%%%%%%%%%%%%%%%%%
\vspace{-0.1cm}
\section{Model setting}
\vspace{-0.3cm}
%%%%%%%%%%%%%%%%%%%%%%%

Consider a continuous-time linear, time invariant, single input single output system,
\begin{equation}\label{eq-System}
	\ \quad x(t) = G_0(p)u(t-\tau_0) =\frac{B(p)}{A(p)}u(t-\tau_0),
\end{equation}
with
\begin{align*}
	B(p) &= b_{0}p^{m}+b_{1}p^{m-1}+...+b_{m}, \\
	A(p) &= p^{n}+a_{1}p^{n-1}+...+a_{n}, \ \ \  n\geq m
\end{align*}
where $\tau_0$ is the time delay, $u(t),\ x(t)$ are the input and deterministic output of the system respectively and $ p $ is the differential operator, i.e. $ p^{(i)}x(t)= d^{i}x(t)/dt^i$. In addition, the following assumptions are made:
\begin{assumption}
Polynomials $B(p)$ and $A(p)$ are coprime.
\end{assumption}
\vspace{-0.4cm}
\begin{assumption}
The system is asymptotically stable.
\end{assumption}
\vspace{-0.4cm}
\begin{assumption} \label{Assum_HFG}
The high frequency gain of $G_0(p)$ is 0, i.e. $G_0(p)$ is strictly proper.
\end{assumption}
\vspace{-0.4cm}

The deterministic output $x(t_k)$ is measured as $y({t_k})$ in the presence of noise, i.e.
\begin{equation}\label{eq-data}
 \ y(t_{k}) = x(t_{k}) + e(t_{k}).
\end{equation}

Furthermore, we consider the sampling time of the input $u(t_k)$ and output data $y(t_k)$ as either regular or irregular. The time-varying sampling interval is denoted as,
\begin{equation}
h_k = t_{k+1}-t_k, \ \ \ \ \ \ \ k = 1,2, ..., N-1,
\end{equation}
where $N$ is the length of the data. 

The objective of a CT system identification problem is to estimate the time delay, $\tau_0$, and the parameters $ a_{1}, a_{2}, ..., a_{n}, b_{0}, b_{1}, ..., b_{m} $ of the CT model in (\ref{eq-System}), using the measured input and output data, $ u(t_{k})$ and $y(t_{k})^{N}_{k=1} $ respectively.
 
%%%%%%%%%%%%%%%%%%%%%%%
\section{Parameter and time delay identification of Continuous-time models with SRIVC and filtering}
%%%%%%%%%%%%%%%%%%%%%%%

There exists a large number of continuous time identification methods see, for example, \cite{Garnier2013} \cite{Young1981}. In this paper, we consider the Simplified Refined Instrumental Variable Continuous Time method (SRIVC), which is developed in the literature by Young and Garnier \cite{Garnier2008}\cite{Garnier2013}\cite{Young1981}. We begin with a basic description of the SRIVC method.

\subsection{Traditional SRIVC method}

To simplify the description we first assume the time delay $\tau_0$ is known. The SRIVC method can then be summarized as follows.

From (\ref{eq-System}) and (\ref{eq-data}),
\begin{equation} \label{eq-sys}
	\frac{A(p)}{A(p)}y(t_{k}) = \frac{B(p)}{A(p)}u(t_{k}-\tau_0) + e(t_{k}),
\end{equation}
which becomes,
\begin{equation} \label{eq-sys-sim}
	A(p)y_{A}(t_{k}) = B(p)u_{A}(t_{k}-\tau_0) + e(t_{k}),
\end{equation}
with
\begin{equation} \label{eq-denote}
\begin{aligned}
& y_{A}(t_{k}) = \dfrac{1}{A(p)}y(t_{k}),\ u_{A}(t_{k}-\tau_0) = \dfrac{1}{A(p)}u(t_{k}-\tau_0).
\end{aligned}
\end{equation}

From (\ref{eq-sys-sim}), a linear regression model can be formed as 
\begin{equation} \label{eq-linear}
	y_{A}^{(n)}(t_{k}) = \varphi_{N}(t_{k},\tau_0)\theta + e(t_{k})
\end{equation}
where 
\begin{equation} \label{eq-phiN}
\begin{aligned}
& \varphi_{N}(t_{k}) = [-y_{A}^{(n-1)}(t_{k}), \ -y_{A}^{(n-2)}(t_{k}), \ ...\ , -y_{A}(t_{k}),
\\
& \ \ \ \ \ \ \ \ \ \ \ \ u_{A}^{(m)}(t_{k}-\tau_0),\  u_{A}^{(m-1)}(t_{k}-\tau_0)\ ...\  u_{A}(t_{k}-\tau_0)], \\
& \theta = [a_{1}, a_{2}, ..., a_{n}, b_{0}, b_{1}, ..., b_{m}]^{T}, \\
& y_{A}^{(i)}(t_{k}) = \dfrac{d^iy_A(t_k)}{dt^i},\ u_{A}^{(i)}(t_{k}-\tau_0) = \dfrac{d^iu_A(t_k-\tau_0)}{dt^i}.
 \end{aligned}
\end{equation}

To solve the linear regression problem in (\ref{eq-linear}), the SRIVC algorithm uses an instrumental variable (IV) method. The instruments in the SRIVC method are chosen as the estimated noise free outputs, i.e. the regressor becomes,
\begin{equation} \label{eq-psiN}
\begin{aligned}
& \hat{\varphi}_{N}(t_{k}) = [-\hat{x}_{A}^{(n-1)}(t_{k}), \ -\hat{x}_{A}^{(n-2)}(t_{k}), \ ...\ , -\hat{x}_{A}(t_{k}),\\
& \ \ \ \ \ \ \ \ \ \ u_{A}^{(m)}(t_{k}-\tau_0),\  u_{A}^{(m-1)}(t_{k}-\tau_0)\ ...\  u_{A}(t_{k}-\tau_0)]
\end{aligned}
\end{equation}
where 
\begin{equation} \label{x_denote}
	\hat{x}_{A}(t_{k}) = \frac{1}{\hat{A}(p)}\hat{x}(t_{k});\ \hat{x}(t_{k})= \frac{\hat{B}(p)}{\hat{A}(p)}u(t_{k}-\tau_0),
\end{equation}
with $\hat{B}(p)$ and $\hat{A}(p)$ estimates of $B(p)$ and $A(p)$.

\subsection{Implementation of the CT filtering operation for irregular sampled data and arbitrary time delay}

The SRIVC method can be used to estimate system parameters from regular sampled data as well as from irregular sampled data. However, one difficulty in implementing the SRIVC method for irregular sampled data in the presence of an arbitrary time delay is the CT filtering operation, e.g., 
\begin{equation} \label{eq_CTfilter}
 u_{A}(t_{k}-\tau) = \dfrac{1}{A(p)}u(t_{k}-\tau),
\end{equation}
when $\tau$ is an arbitrary time delay.

There are two reasons for this difficulty:
\begin{enumerate}
\item
$u(t_k-\tau)$ is not available from the measured data,
\item
The digital simulation of $u(t_k-\tau)$ is generally performed in state-space form hence, an equivalent discrete-time state space representation of the CT state space model needs to be computed. For irregular data, as the sampling interval $h(t_k)$ is time-varying, the computational load using standard methods, e.g. \textit{expm} in Matlab, to compute the transformation matrices will be large.
\end{enumerate}

As suggested in \cite{Chen2015}, the two problems mentioned above can be solved as follows,

\begin{enumerate}
\item
$u(t_k-\tau)$ can be constructed from the neighbouring data based on the inter-sample behaviour, e.g. zero-order-hold (ZOH), first-order-hold (FOH), etc.
\item
To reduce the computational load, the time-varying sampling interval is divided into two intervals: the first interval is a multiple of a constant sampling period, the second interval is the residual. The discretization matrix of the first interval is pre-computed using a standard method and stored in an array. The discretization matrix for the second interval is computed by a fast approximation, e.g. the 4th order Runge-Kutta (RK4) method. The final discretization matrix of the sampling interval $h_k$ is a product of the two matrices.

Specifically, convert the transfer function (TF) model in (\ref{eq_CTfilter}) to the equivalent CT state-space model,
\begin{equation}
\begin{cases}
    \dot{\textbf{z}}(t-\tau) = \textbf{F}_c\textbf{z}(t-\tau)+\textbf{G}_cu(t-\tau) \\
    u_A(t-\tau) = \textbf{H}_c\textbf{z}(t-\tau).
\end{cases}
\end{equation}
The equivalent DT state-space model will be,
\begin{equation} \label{ss_dt_der}
\begin{cases}
    \textbf{z}(t_{k+1}-\tau) = \textbf{F}(h_k)\textbf{z}(t_k-\tau)+\textbf{M}(t_k) \\
    u_A(t_k-\tau) = \textbf{H}_c\textbf{z}(t_k-\tau)
\end{cases}
\end{equation}
where $\textbf{F}(h_k)$, $\textbf{M}(t_k)$ are computed as follows,
\begin{enumerate}
\item
Denote $t_{k,1}, t_{k,n_k}$ ($n_k \geq 2$) as $t_k-\tau$, $t_{k+1}-\tau$ respectively; $t_{k,i}, i = 2, ..., n_k-1$ as the transition time-instants of $u(t)$ between $t_k-\tau$ and $t_{k+1}-\tau$; and $h_{k,i} = t_{k,i+1}-t_{k,i}$.
\item
Let $\Delta$ be a constant sampling period; $m_{k,i}$ be a positive integer; $\delta_{k,i} \geq 0$ such that:
\begin{center}
$h_{k,i}=m_{k,i}\Delta+\delta_{k,i}<(m_{k,i}+1)\Delta$.
\end{center}
\item
Compute $\textbf{F}(m\Delta), \textbf{G}(m\Delta)$.
\item
Compute $\textbf{F}(\delta_{k,i}), \textbf{G}(\delta_{k,i})$ using, e.g. RK4.
\item
Finally, compute $\textbf{F}(h_k)$, $\textbf{M}(t_k)$ as follows,
\begin{equation} \nonumber
\begin{aligned}
& \textbf{F}(h_k) = e^{\textbf{F}_ch_k}, \\
& \textbf{M}(t_k) = \sum\limits_{i=1}^{n_k-1} \prod\limits_{j=i+1}^{n_k-1} \textbf{F}(h_{k,j})\textbf{G}(h_{k,i}) u(t_{k,i}), \\
\end{aligned}
\end{equation}
with
\begin{equation} \nonumber
\begin{aligned}
& \textbf{F}(h_{k,i}) = \textbf{F}(m_{k,i}\Delta)\textbf{F}(\delta_{k,i}), \\
& \textbf{G}(h_{k,i}) = \textbf{F}(m_{k,i}\Delta)\textbf{G}(\delta_{k,i}) + \textbf{G}(m_{k,i}\Delta).\\
\end{aligned}
\end{equation}
\end{enumerate}
\end{enumerate}

Now that we can generate filtered signals of the irregular sampled data, we describe the traditional SRIVC algorithm \cite{Garnier2008}\cite{Garnier2013}\cite{Young1981} in Algorithm I.\\
{\begin{center}
\textsc{ALGORITHM I} 
\end{center}}
\vspace{0.4cm}
\begin{mdframed}
\textbf{Step 1. Initialization} 
\begin{enumerate}
\item
Create the stable state variable filter (SVF),
\begin{equation} \label{eq-svf}
	F(p) = \dfrac{1}{(p+\omega_c^{SVF})^n},
\end{equation}
where $\omega_c^{SVF}$ is chosen to be larger than or equal to the bandwidth of the system.
\item
Filter $y(t_{k})$ and $u(t_{k}-\tau_0)$ via the SVF to generate derivatives of the signals, i.e. $y_{A}^{(n)}(t_{k}), y_{A}^{(n-1)}(t_{k}), ...,$ $y_{A}(t_{k})$ and $u_{A}^{(m)}(t_{k}-\tau_0), u_{A}^{(m-1)}(t_{k}-\tau_0), ..., u_{A}(t_{k}-\tau_0)$ by using the method described in \cite{Garnier2008}.
\item
Use the least squares method to estimate the initial parameters,

\begin{equation}
\hat{\theta}_{0} = [\Phi_{N}\Phi_{N}^T]^{-1}\Phi_{N}Y_{N}.
\end{equation}
with
\begin{equation}
\begin{aligned}
& \Phi_{N} = [\varphi_{N}^T(t_{1})\ ...\ \varphi_{N}^T(t_{N})],
\\
& Y_{N} = [y_{A}^{(n)}(t_{1})\ y_{A}^{(n)}(t_{2})\ ...\ y_{A}^{(n)}(t_{N})]^T,
\end{aligned}
\end{equation}
where $\varphi_N$ is defined in (\ref{eq-phiN}).

\end{enumerate}
\textbf{Step 2. Iterative estimation}

\textit{for} j = 1:\textit{convergence}\footnote{Convergence requires the relative error between the estimated parameters of two consecutive iterations to be less than an $\epsilon>0$ where $\epsilon$ is a small value chosen to achieve the desired accuracy.}
\begin{enumerate}
\item
If the estimated TF model is unstable, reflect the right half plane zeros of $\hat{A}(p,\hat{\theta}^{j-1})$ into the left half plane and construct the new estimate. 
\item
Generate the instrumental variable,
\begin{equation}
 \hat{x}(t_{k}) = \dfrac{\hat{B}(p,\hat{\theta}^{j-1})}{\hat{A}(p,\hat{\theta}^{j-1})}u(t_{k}-\tau_0).
\end{equation}
\item
Filter the input $u(t_{k}-\tau_0)$, output $y(t_{k})$ and the instrumental variable $\hat{x}(t_{k})$ by the continuous time filter
\begin{equation}
 F_{c}(p) = \dfrac{1}{\hat{A}(p,\hat{\theta}^{j-1})},
\end{equation}
to generate the derivatives of these signals.
\item
Using the prefiltered data, generate an estimate $\hat{\theta}^{j}$ using the IV method,
\begin{equation} \label{eq-iv}
\hat{\theta}^{j} = [\Psi_{N}\Phi_{N}^T]^{-1}\Psi_{N}Y_{N}.
\end{equation}
where $\Psi_{N}$ is the IV matrix generated by the instrumental variables $\hat{x}(t_{k})$ and $\Phi_{N}$ is the regression matrix,
\begin{equation}
\begin{aligned}
& \Phi_{N} = [\varphi_{N}^T(t_{1})\ ...\ \varphi_{N}^T(t_{N})],
\\
& \Psi_{N} = [\hat{\varphi}_{N}^T(t_{1})\ ...\ \hat{\varphi}_{N}^T(t_{N})],
\end{aligned}
\end{equation}
where $\hat{\varphi}_N$ is defined in (\ref{eq-psiN}).
\end{enumerate} 
\textit{end}
\end{mdframed}
\vspace{0.2cm}

\subsection{SRIVC-based time delay estimation with filtering}
A recently developed method \cite{Chen2015} to estimate both the time delay and parameters for the CT model (\ref{eq-data}) considers it as a separable nonlinear least squares problem. The SRIVC algorithm and the Gauss-Newton method are used in \cite{Chen2015} to estimate the system parameters and the time delay respectively. In this problem the cost function for the time delay estimation has multiple minima \cite{Chen2015}, hence a low-pass filter is utilized to extend the global convergence region \cite{Chen2015} \cite{Diego2016} \cite{Ferretti:96}.

Let $L(p)$ be a CT low-pass filter with cut off frequency $\omega_c^L$. Applying the filter $L(p)$ to the linear regression (\ref{eq-linear}), we obtain
\begin{equation} \label{eq-linearFilter}
	\bar{y}_{A}^{(n)}(\rho) = \bar{\varphi}_{N}(\rho)\theta + \bar{\epsilon}(\rho),
\end{equation}
where $\rho = [\theta^T, \tau]$ and $\bar{.}$ represents the signal $y,u$ or $e$ filtered by $L(p)$.

When the estimated parameter $\hat{\theta}$ is a function of $\tau$, i.e.
\begin{equation}
\hat{\theta}(\tau) = [\bar{\Psi}_{N}(\rho)\bar{\Phi}_{N}^T(\rho)]^{-1}\bar{\Psi}_{N}(\rho)\bar{Y}_{N},
\end{equation}
then the time delay $\tau$ can be estimated as \cite{Chen2015}\cite{Ngia2001},
\begin{equation}
\hat{\tau} = \underset{\tau}{\text{arg min}} \tilde{\bar{J_N}}(\tau) = \underset{\tau}{\text{arg min}} \bar{J_N} (\rho) \vert_{\theta=\hat{\theta}(\tau)}, 
\end{equation}
with
\begin{equation}
\bar{J}_N(\rho) = \dfrac{1}{2(N-s+1)}\bar{\epsilon}^T(\rho)\bar{\epsilon}(\rho),
\end{equation}
where $s$ is chosen to guarantee that at sample $s$, $t_s \geq \tau$.

By using the Gauss-Newton method, the system parameters and the time delay can be iteratively estimated,
\begin{equation} \label{Gauss_est}
\begin{aligned}
\hat{\tau}^{j+1} &= \hat{\tau}^j-\mu^j\ [\nabla^2\tilde{\bar{J_N}}(\hat{\tau}^j)] ^{-1} \nabla \tilde{\bar{J_N}}(\hat{\tau}^j) \\
\hat{\theta}^{j+1} &= [\bar{\Psi}_{N}(\hat{\theta}^{j},\hat{\tau}^{j+1})\bar{\Phi}_{N}^T(\hat{\theta}^{j},\hat{\tau}^{j+1})]^{-1} \bar{\Psi}_{N}(\hat{\theta}^{j},\hat{\tau}^{j+1})\bar{Y}_{N}(\hat{\theta}^{j}),
\end{aligned}
\end{equation}
where $\mu^j$ is the step size and
\begin{equation} \label{Gauss_def}
\begin{aligned}
& \nabla\tilde{\bar{J_N}}(\hat{\tau}^j) = \bar{\epsilon}^T_{\tau}\bar{\epsilon}(\rho),
\\
& \nabla^2\tilde{\bar{J_N}}(\hat{\tau}^j) = \bar{\epsilon}^T_{\tau}\bar{\epsilon}_{\tau}-\bar{\epsilon}^T_{\tau}\bar{\epsilon_{\theta}}(\bar{\epsilon}^T_{\theta}\bar{\epsilon_{\theta}})^{-1}\bar{\epsilon}^T_{\theta}\bar{\epsilon_{\tau}},
\\
& \epsilon_{\tau} = \dfrac{\partial \epsilon(\rho)}{\partial\tau} \vert_{\rho=\hat{\rho}^j} = pG(p,\hat{\theta}^j)u(\hat{\tau}^{j}),
\\
& \epsilon_{\theta} = \dfrac{\partial \epsilon(\rho)}{\partial\theta} \vert_{\rho=\hat{\rho}^j} = -\bar{\Psi}_{N}^T(\hat{\theta}^{j},\hat{\tau}^{j}).
\end{aligned}
\end{equation}

The SRIVC-based time delay estimation with filtering \cite{Chen2015} is summarized in Algorithm II. \\
{\begin{center}
\textsc{ALGORITHM II} 
\end{center}}
\vspace{0.4cm}
\begin{mdframed}
\textbf{Step 1. Initialization} 
\begin{enumerate}
\item
Set the initial value $\hat{\tau}^0$, the boundaries\footnote{$\tau_{min}$ and $\tau_{max}$ are boundaries for the delay, which are known a priori, i.e. $\tau_{min} \leq \tau_0 \leq \tau_{max}$} $\tau_{min}, \tau_{max}$, $\Delta \tau_{min}$\footnote{$\Delta \tau_{min}$ is the limit value of the step size of the time delay estimate using Gauss-Newton method.}, the cut-off frequency $\omega_c^L$ of $L(p)$, the cut-off frequency of $F(p)$, $\omega_c^{SVF}$, and a small positive $\epsilon$ to indicate convergence.
\item
Based on the initial value $\hat{\tau}^0$, use the SRIVC method to compute $\hat{\theta}^0$.
\item
Compute $\nabla^2\tilde{\bar{J_N}}(\hat{\tau}^0)$ and $\nabla \tilde{\bar{J_N}}(\hat{\tau}^0)$ from (\ref{Gauss_def}).
\end{enumerate}
\textbf{Step 2. Iterative estimation}

\textit{for} j = 1:\textit{convergence}
\begin{enumerate}
\item
Compute $\Delta\hat{\tau}^j$ using the filtered input/output data, equation (\ref{Gauss_def}) and \begin{equation}
\Delta\hat{\tau}^j = - [\nabla^2\tilde{\bar{J_N}}(\hat{\tau}^j)] ^{-1} \nabla \tilde{\bar{J_N}}(\hat{\tau}^j).
\end{equation}
\item

\begin{enumerate}
\item
Compute $\hat{\tau}^{j+1} = \hat{\tau}^j + \Delta\hat{\tau}^j$
\\
If $\hat{\tau}^{j+1} \not\in [\tau_{min},\tau_{max}]$, let $\Delta\hat{\tau}^j = \Delta\hat{\tau}^j/2$ and repeat this step.\\
If $\vert \Delta\hat{\tau}^j \vert \leq \Delta \tau_{min}$, break.
\item
Estimate $\hat{\theta}^{j+1}$ using the SRIVC method and time delay $\hat{\tau}^{j+1}$ from the filtered input/output data.
\item
Compute $\tilde{\bar{J_N}}(\hat{\tau}^{j+1})$. If $\tilde{\bar{J_N}}(\hat{\tau}^{j+1}) \geq \tilde{\bar{J_N}}(\hat{\tau}^{j})$, let $\Delta\hat{\tau}^j = \Delta\hat{\tau}^j/2$ and return to (a).
\end{enumerate}
\item
If $\tilde{\bar{J_N}}(\hat{\tau}^{j})-\tilde{\bar{J_N}}(\hat{\tau}^{j+1}) \geq \epsilon$, go to Step 1, else break.
\end{enumerate}
\ \ \ \   \textit{end}
\\
\textbf{Step 3. Refined parameter estimation}

Repeat Step 2 with low-pass filter $L(p) = 1$.

\end{mdframed}
\begin{remark} \label{rm_coFreq}
The smaller the cut-off frequency $\omega_c^L$, the larger the global convergence region. However, when $\omega_c^L$ is chosen too small, information from the data is lost, which affects the accuracy of the parameter estimation \cite{Chen2015}. The suggestion in \cite{Chen2015} is, 
\begin{equation} \label{eq_wc2}
\omega_c^L \geq \dfrac{1}{10}\text{bw},
\end{equation}
where $\text{bw}$ is the bandwidth of the system $G_0(p)$ (rad/sec).
\end{remark}

%%%%%%%%%%%%%%%%%%%%%%%
\section{SRIVC-based time delay estimation with multiple filtering and redundancy method}
%%%%%%%%%%%%%%%%%%%%%%%

In this section, we propose a method that improves convergence to the global minimum of the time delay optimization problem using the useful redundancy technique \cite{Alamir2008}\cite{Alamir2009}. The useful redundancy technique was developed to avoid local minima when solving a non-linear inverse problem. As stated in the previous section, the main problem in time delay estimation is that the cost function $J_N(\rho)$ possesses many local minima \cite{Chen2015}. The filtering operation \cite{Ferretti:96} described in Section 3.3 can be employed in order to increase the global convergence region of the cost function associated with the time delay estimation. However, when the data is very noisy or the initial value, $\hat{\tau}^0$, is located far from the global minimum, using only one filter does not guarantee that the solution of the optimization problem converges to the global minimum. 

To demonstrate how useful redundancy can be utilized in our estimator, we first describe the useful redundancy technique.

\subsection{The useful redundancy technique}

We define the useful redundancy technique by quoting directly from \cite{Alamir2009}.

\textbf{Definition 1} \cite{Alamir2009}. Consider an  optimization problem,
\begin{equation} \nonumber
\underset{\rho}{\text{min}} \ \ J_0(\rho).
\end{equation}
Then it is called $\textit{M-safely redundant}$ if and only if the following conditions hold:
\begin{enumerate}
\item
There exists a finite $M$ cost functions $J_i$ sharing the same global minimum $\rho^* \in \mathbb{R}^n$.
\item
There exists a solver (or an iterative scheme) $\ell$ and a finite number of iterations $r^* \in \mathbb{N}$ such that for some $\gamma \in [0,1)$ and all $\rho \in \mathbb{P}$ the following inequality holds,
\begin{equation}
\Delta^{\gamma}_N(\rho) = \underset{i \in \lbrace 0,...,M \rbrace}{\text{min}} [J_0(\ell^{(r^*)}(\rho,J_i))-\gamma J_0(\rho)]\ \leq 0
\end{equation}
where $\ell^{(r^*)}(\rho,J_i)$ is the candidate solution obtained after $r^*$ iterations of $\ell$ using the cost function $J_i$ and starting from the initial guess $\rho$. $ \hfill\square$ 
\end{enumerate}
\vspace{0.1cm}

The solver path $(\rho, J_i)$ is defined as the sequence of iterates $\ell^{(j)}(\rho,J_i)$ for the solver $\ell$ when the cost function $J_i$ starts from an initial guess $\rho$. Condition 2 means that for any initial $\rho$, there always exists a solver path $\ell$ that corresponds to a decrease in the original cost $J_0$ after at most $r^*$ iterations. It is proven in \cite{Alamir2009} that if an optimization problem is \textit{M-safely redundant} following from Definition 1, that convergence to a global minimum is guaranteed.

An algorithm that describes the useful redundancy technique is given in Algorithm III. 
\vspace{0.2cm}
{\begin{center}
\textsc{ALGORITHM III} 
\end{center}}
\vspace{0.4cm}
\begin{mdframed}
\textbf{Step 1. Initialization} 
Choose an initial value $\rho^{(0)}$. \\
\textbf{Step 2. Iterative estimation} \\
\textit{for} \textit{$k$=0:converge}
\begin{enumerate}
\item
Set $i=0$ and $failure = true$
\item
While \textit{$failure$}  do
\begin{enumerate}
\item
Use the solver path $(\rho, J_i)$ starting from $\rho^{(k)}$, find the corresponding minima $\xi^{(k,i)}$.
\item
Compute \\
$failure = \Big(J_0(\xi^{(k,i)}) > (1-\epsilon)J_0(\rho^{(k)}) \Big)$ \footnote{$\epsilon$ is a predefined small value that is chosen based on the desired accuracy.}
\item
If $failure$ then $i = (i+1)\ \text{mod}\ M$, \\
Else set $k = k+1, \rho^{(k)} = \xi^{(k,i)}$.
\end{enumerate}
\text{End while}
\end{enumerate}
\textit{end}
\end{mdframed}

\vspace{0.1cm}
\subsection{Theoretical results related to the time delay cost function}
\vspace{-0.1cm}

To construct an \textit{M-safely redundant} optimization problem for the estimation of the time delay, we need a cost function $J_0$ and multiple solver paths that satisfy the two conditions in Definition 1. Here, the solver paths are generated by filtering the time delay estimation error using a set of low-pass filters with different cut-off frequencies. The cost function $J_0$ is formulated from these filtered cost functions such that there always exists a solver path that corresponds to a decrease in $J_0$ after a finite number of iterations. Next we describe the set of filters and the cost function $J_0$ required to satisfy the conditions in Definition 1.

Consider a continuous-time, linear, time-invariant SISO system described by
\begin{equation} \label{md_setting}
y(t) = G_0(p)u(t-\tau_0) + e(t).
\end{equation}
We make a further assumption on the noise $e(t)$. 
\begin{assumption} \label{Assum_noise}
$e(t)$ is white random process uncorrelated with $u$ having intensity $\lambda$.
\end{assumption}

For an estimate $G(p,\theta)$ and $\tau$, the estimation error $\epsilon(t,\theta,\tau)$ can be computed as,
\begin{equation}
\epsilon(t,\theta,\tau) = y(t) - G(p,\theta)u(t-\tau).
\end{equation}

As mentioned in Section 3, the delay can be estimated by minimizing the cost function $J(\theta,\tau)$, where $J(\theta,\tau) = \int_{-\infty}^{\infty} \! \epsilon(t,\theta,\tau)^2 \, \mathrm{dt} $. If the estimation error is filtered by the low-pass filter $L(p)$, then an estimate of $\theta$ and $\tau$ can be computed by,
\begin{equation} \label{ls_problem}
(\hat{\theta},\hat{\tau}) = \text{arg min}\ \bar{J}(\theta,\tau),
\end{equation}
where
\begin{equation} \label{J_filter}
\begin{aligned}
\bar{J}(\theta,\tau) &= \int_{-\infty}^{\infty} \! \bar{\epsilon}(t,\theta,\tau)^2 \, \mathrm{dt} \ = \int_{-\infty}^{\infty} \! \Big\lbrace L(p)\epsilon(t,\theta,\tau) \Big\rbrace^2 \, \mathrm{dt} \\
&= \int_{-\infty}^{\infty} \! \Big\lbrace L(p)[y(t)-G(p,\theta)u(t-\tau)] \Big\rbrace^2 \, \mathrm{dt} \\
\end{aligned}
\end{equation}
which by Parseval's theorem is equivalent to,
\begin{equation} \label{J_filter}
\begin{aligned}
\bar{J}(\theta,\tau) &= \dfrac{1}{2\pi}\int_{-\infty}^{\infty} \! \Big[|G_0(j\omega)e^{-j\omega\tau_0}-G(j\omega,\theta)e^{-j\omega\tau}|^2 \\
& \ \ \ \ \ \ \ \ \ \ \ \ \ \ \ \ \ \ \times\Psi_u(\omega) +\Psi_v(\omega)\Big] \big\vert L(j\omega)\big\vert^2 \, \mathrm{d}\omega, \\
\end{aligned}
\end{equation}
with $\Psi_u(\omega)$ the spectral density of $u(t)$ and the spectral density $\Psi_v(\omega) = \lambda$ follows from Assumption \ref{Assum_noise}. If the transfer function $G_0(p)$ is known; the input signal white noise, i.e. $\Psi_u(\omega) = 1$; and $\delta\tau = \tau-\tau_0$, the cost function $\bar{J}(\theta,\tau)$ can be written as,
\begin{equation} \label{J_filter2}
\begin{aligned}
& \bar{J}(\delta\tau) = \dfrac{1}{2\pi}\int_{-\infty}^{\infty} \! \Big[ \vert 1-e^{-j\omega\delta\tau} \vert^2 \big\vert G_0(j\omega)e^{-j\omega\tau_0} \vert^2\\
& \ \ \ \ \ \ \ \ \ \ \ \ \ \ \ \ \ \ \ \ \ \ \ \ \ \ \times \Psi_u(\omega) +\lambda\Big] \big\vert L(j\omega)\big\vert^2 \, \mathrm{d}\omega \\
& \ \ \ \ \ \ \ = \dfrac{1}{\pi}\int_{-\infty}^{\infty} \! \Big[ (1-\cos(\omega\delta\tau)) \big\vert G_0(j\omega)\vert^2 +\dfrac{\lambda}{2}\Big] \\
& \ \ \ \ \ \ \ \ \ \ \ \ \ \ \ \ \ \ \ \ \ \ \ \ \ \ \times \big\vert L(j\omega)\big\vert^2 \, \mathrm{d}\omega. 
\end{aligned}
\end{equation}
Recall that we are concerned with how to choose a set of filters and the cost function $J_0$ such that the two conditions of an M-safely redundant problem are satisfied. To satisfy the first condition we need to ensure all the cost functions share the same global minimum. The second condition requires the filter set to be constructed such that for any initial value of $\delta\tau$, there always exists a cost function whose solver path corresponds to a decrease in the cost function $J_0$ after a finite number of iterations. Note that these conditions can be checked if we know the locations of the global minimum and the extrema, i.e. the minima and maxima of $\bar{J}(\delta\tau)$. 

First we establish a result for the location of the global minimum of the (non)filtered delay cost function. 
  
\begin{theorem} \label{theorem_gbm}
Consider the system $G_0(p)$ as described in (\ref{md_setting}). When $\lambda = 0$, $\forall \delta\tau \in \mathbb{R}$, $\forall$ low-pass filters $L(p) \neq 0$, such that $\bar{J}(\delta\tau) \geq \bar{J}(0)$, the equality $\bar{J}(\delta\tau) = \bar{J}(0)$ occurs if and only if $\delta\tau = 0$.
\end{theorem}
\textbf{Proof.} The proof of Theorem \ref{theorem_gbm} is provided in Appendix A.1.

Theorem \ref{theorem_minloc} provides a result for the locations of the extrema of the filtered cost function $\bar{J}(\delta\tau)$. 
\begin{theorem} \label{theorem_minloc}
If $L(p)$ is selected such that $L(p)G_0(p)$ is an ideal low-pass filter with cut-off frequency $\omega_c$ rad/sec, then the locations of the $i^{th}$ positive extrema of the time delay cost function $\bar{J}(\delta\tau)$ can be approximated by,
\begin{equation} \label{eq_minloc}
\begin{aligned}
& \tilde{\delta}\tau_i \simeq \dfrac{2i+1}{4}T_c - \dfrac{T_c}{(2i+1)\pi^2},
\end{aligned}
\end{equation}
where $T_c = 2\pi/\omega_c$ and $i \geq 1$. When $i$ is even, the corresponding extremum is a minimum and when $i$ is odd, it is a maximum. 
\end{theorem}
\textbf{Proof.} The proof of Theorem \ref{theorem_minloc} is provided in Appendix A.2.

\begin{remark} \label{rm_gconverg}
From Theorem \ref{theorem_minloc}, the locations of the extrema for the filtered cost function are known, hence the global convergence region is, 
\begin{equation} \label{eq_gbcoverg}
|\delta\tau| < \dfrac{3\pi}{2\omega_c}-\dfrac{2}{3\pi \omega_c},
\end{equation}
corresponding to the distance from the maxima closest to the global minimum. 
\end{remark}

\begin{figure}[H]
\begin{center}
\includegraphics[width=73mm,height=51mm]{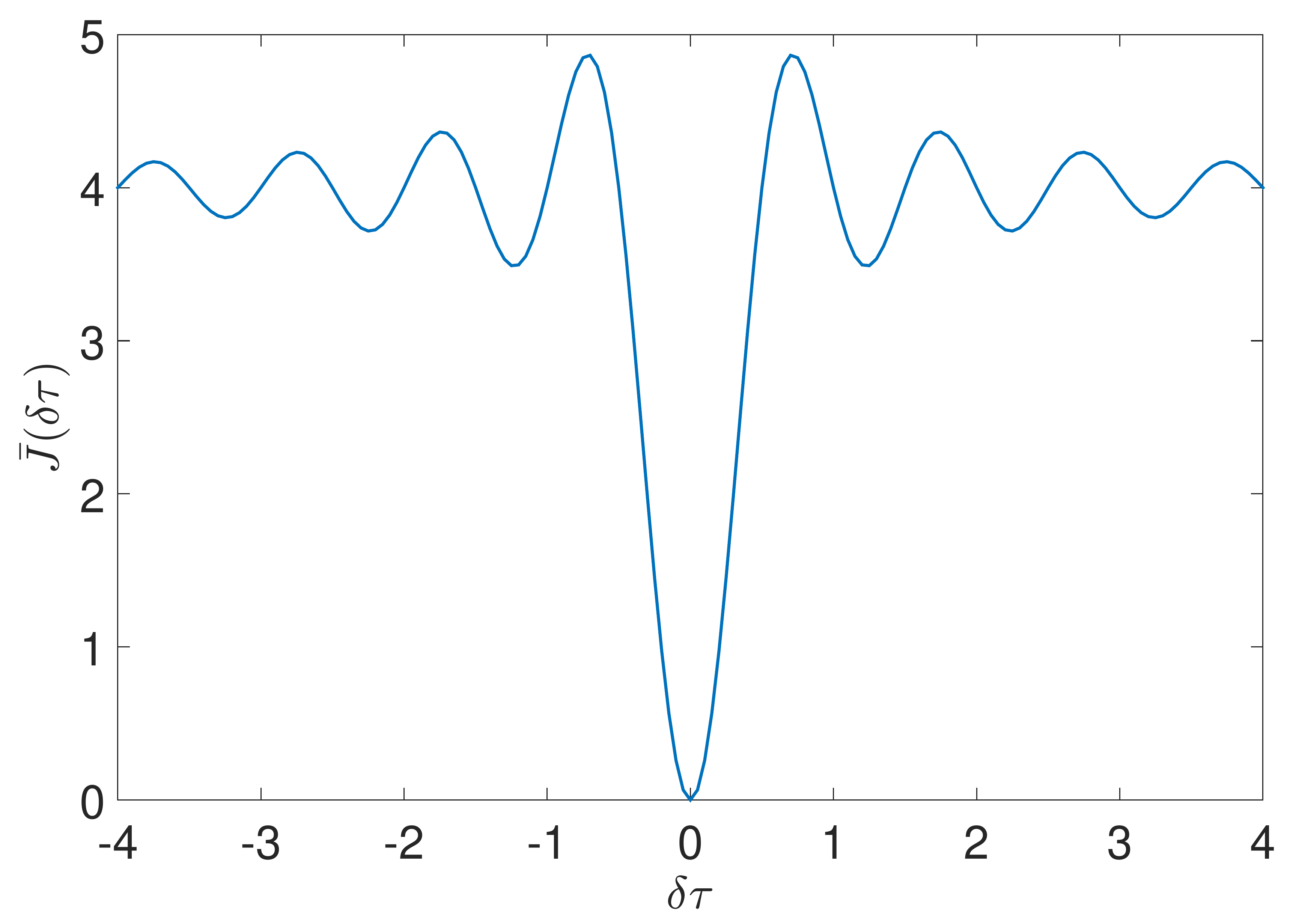}
\caption{Plot of $ \bar{J}(\delta\tau) $ when $\omega_c=2\pi$} 
\label{Simple_localMinima} 
\end{center}
\end{figure}

The result of Remark \ref{rm_gconverg} can be seen in Fig. \ref{Simple_localMinima}, which provides the plot of $\bar{J}(\delta\tau)$ when $\omega_c=2\pi$ (the constant $\lambda$ of the integral in (\ref{J_filter2}) is set to 0 here for simplicity). 

\subsection{Time delay estimation with the useful redundancy technique}

Now we can specify a filter set that satisfies the two conditions in Definition 1.

\textit{Condition 1:}
From Theorem \ref{theorem_gbm}, the first condition is satisfied for any set of low-pass filters as the global minimum of all the filtered time delay cost functions occurs when $\delta\tau = 0$.

\textit{Condition 2:}
The second condition depends critically on the choice of filters and cost function $J_0$. The filter set needs to be chosen such that for any initial value of $\delta\tau$, there always exists a cost function whose solver path corresponds to a decrease in the cost function $J_0$ after a finite number of iterations. 

To verify that condition 2 is satisfied, consider the necessary and sufficient conditions. To satisfy the necessary condition, we need to show that for any initial value $\delta\tau_0$, there always exists a path that has a minima location, $\delta\tau_1$, closer to the global minimum w.r.t. $\delta\tau_0$. For the sufficient condition, we need to establish that $J_0(\delta\tau_1)<J_0(\delta\tau_0)$.

First, we check the necessary condition. A simple choice is a filter set $L_k(p), k = \overline{1,n_f}$ where the cut-off frequencies, $\omega_{c,k}$, of $L_k(p)G_0(p)$ span from $1/\beta$ to 1 of the system bandwidth ($10 \geq \beta > 1$)\footnote{From Remark \ref{rm_coFreq}, $\beta$ is chosen $\leq10$.}. The cut-off frequencies are chosen based on the corresponding periods, $T_{c,k}$, being linearly spaced, i.e.,
\begin{equation}
T_{c,k} = \beta T_{bw} + (k-1)\dfrac{(1-\beta)T_{bw}}{n_f-1}, \ \ \ \ k = \overline{1,n_f},
\end{equation}
where $T_{bw} = 2\pi/bw$, $bw$ is the bandwidth of the system $G_0(p)$ (rad/sec). This is due to the local minima of the filtered cost function being linearly dependant on the period of the cut-off frequency (Theorem \ref{theorem_minloc}). 

Finally we need to prove the existence of such a filter set that satisfies the necessary condition of Condition 2 in Definition 1.

\begin{theorem} \label{theorem_nf}
There exists a filter set $L_k(p),\ k=\overline{1,n_f}$ chosen such that, $\forall \ \beta \geq \dfrac{5/4-1/(5\pi^2)}{3/4-1/(3\pi^2)}$, where $L_k(p)G_0(p)$ are ideal low-pass filters chosen with linearly spaced periods such that the cut-off frequencies, $\omega_{c,k}$, spanning from $1/\beta$ to 1 of the system bandwidth, i.e.
\begin{equation}
T_{c,k} = \beta T_{bw} - (k-1)\dfrac{(\beta-1)T_{bw}}{n_f-1}, \ \ \ \ k = \overline{1,n_f},
\end{equation}
where $T_{bw} = 2\pi/bw$ and $T_{c,k} = 2\pi/\omega_{c,k}$. Then,
\begin{equation} \label{eq_t3_cond}
\forall \delta\tau_0 \neq 0, \exists L_q(p): |\xi(L_q(p),\delta\tau_0)| < |\delta\tau_0|,
\end{equation}
where $\xi(L_q(p),\delta\tau_0)$ is the corresponding minimum found using the filtered time delay cost function generated by $L_q(p)$ with the initial delay $\delta\tau_0$. 
\end{theorem}

\textbf{Proof.} The proof of Theorem \ref{theorem_nf} is provided in Appendix A.3.

\begin{remark}
From Theorem \ref{theorem_nf}, we see that with the choice of,
\begin{equation}
\beta \geq \dfrac{5/4-1/(5\pi^2)}{3/4-1/(3\pi^2)},
\end{equation}
there always exists an $n_f$, that can be computed by (\ref{eq_t3_nf}) in Appendix A.3, that ensures the filter set defined in Theorem \ref{theorem_nf} satisfies the necessary part of Definition 1, Condition 2. 
\end{remark}

\begin{remark}
Note that (\ref{eq_t3_nf}) provides a loose lower bound. It is possible to have a smaller value of $n_f$ and still obtain a filter set that satisfies the necessary condition. 
\end{remark}

\begin{remark}
If $G_0(p)$ does not have resonant peaks then it doesn't matter if $G_0(p)$ is known or unknown, as any ideal low pass filter $L_k(p)$ with bandwidth smaller than the bandwidth of $G_0(p)$ will allow $L_k(p)G_0 (p)$ to approximate the desired low pass filter behaviour. If $G_0(p)$ has resonant peaks and the bandwidth of $L_k(p)$ is chosen significantly smaller than that of $G_0(p)$ then $L_k(p)G_0 (p)$ will approximate the desired low pass filter behaviour. 
\end{remark}

\begin{figure}[h]
\begin{center}
\includegraphics[width=75mm,height=52mm]{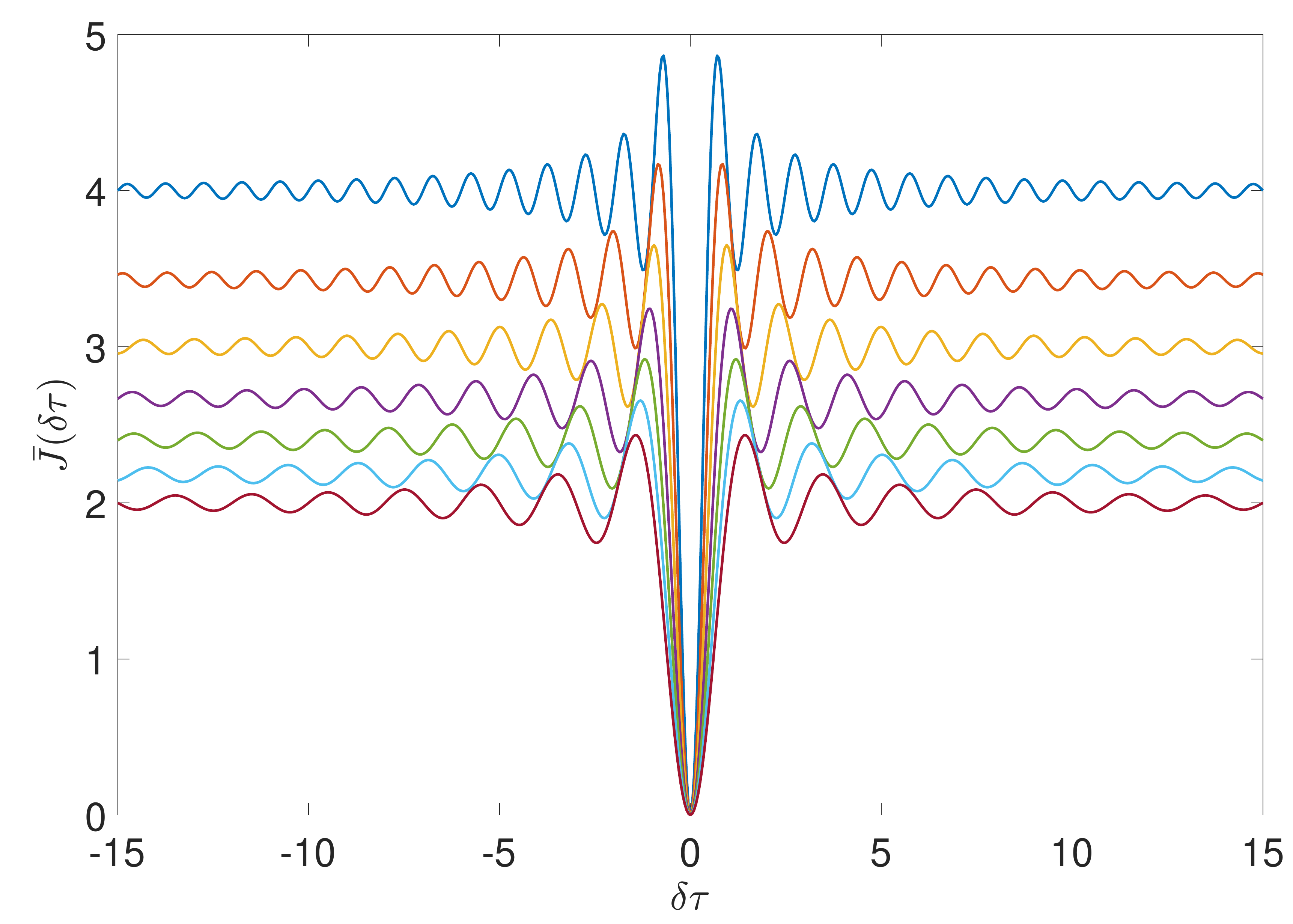}
\caption{Plots of $J(\delta\tau)$ set of filters with linearly spaced periods} 
\label{Simple_setFilter} 
\end{center}
\end{figure}

An empirical method can also be utilized to determine the number of filters. For example, to check if a set of $n_f$ filters with the ratio $\beta$ satisfies the necessary condition, we can plot $n_f$ cost functions that follow from (\ref{J_eq}) in Appendix A.2 and check the minima locations using (\ref{eq_minloc}) to see that for any initial delay, $\delta\tau_0$, there always exists a cost function where the corresponding minima $\delta\tau_i$ is found closer to the global minimum with respect to $\delta\tau_0$. Fig. \ref{Simple_setFilter} shows a set of six cost functions with the ratio $\beta = 2$ for a system bandwidth of $2\pi$ rad/sec, i.e. the cut-off frequencies span from $1/2$ to $1$ (Hz). From the figure, we can see that for any initial value of the time delay less than 15 sec, there always exists a cost function where the minima is found closer to the global minimum with respect to the initial delay. This can be confirmed by computing the minima of all the cost functions using (\ref{eq_minloc}) and observing the path to the global minimum. 

Note that in practice, as $G_0(p)$ is not known a priori, and with the difficulty of implementing an ideal low-pass filter with irregular data, the filters $L_k(p)$ are implemented as Butterworth filters with cut-off frequencies $\omega_{c,k}$. Due to this, the extrema locations are not exactly those specified in (\ref{eq_minloc}). We discuss the effect of this difference by first citing the following theorem from \cite{Ferretti:96}, as a lemma.
\begin{lemma} \label{lemma_Jdelay}
\cite{Ferretti:96} Consider a system with transfer function $F(s) = L(s)G_0(s)e^{-s\tau_0}$, let $u$ be the input signal and $w$ the corresponding output signal of the system. Then we have the following relationship,
\begin{equation} \label{eq_lemm_Jdelay}
\bar{J}(\delta\tau) = 2[\gamma_w(0)-\gamma_w(\delta\tau)] = 2\gamma_w(0)[1-\rho_w(\delta\tau)],
\end{equation}
where $\gamma_w(t)$ is the autocorrelation function of $w(t)$ and $\rho_w(\delta\tau) = \gamma_w(\delta\tau)/\gamma_w(0)$.
\end{lemma}
\textbf{Proof.} See proof in \cite{Ferretti:96}.

With respect to the model setting in (\ref{md_setting}), $\gamma_w(\delta\tau)$ in Lemma \ref{lemma_Jdelay} is the autocorrelation of the output of the system. Hence from Lemma \ref{lemma_Jdelay}, for an input signal of Gaussian distributed white noise, the filtered delay cost function $\bar{J}(\delta\tau)$ will have the form,
\begin{equation} \label{aucorr_form}
\bar{J}(\delta\tau) = \sum\limits_{k=1}^{n_1} a_ke^{-b_k\delta\tau} + \sum\limits_{k=1}^{n_2} h_ke^{-c_k\delta\tau}\cos(d_k\delta\tau+\psi_k),
\end{equation}
where $b_k$ is the absolute value of the $k^{th}$ real pole of $L(s)G_0(s)$; $c_k$ and $d_k$ are the absolute values of the real and imaginary parts of the $k^{th}$ complex pole of $L(s)G_0(s)$; and weighting coefficients $h_k$ are related to the magnitude of the frequency response at the corresponding frequency. Augmenting a low-pass filter $L(s)$ to $G_0(s)$ obviously introduces slow poles into the system and by choosing $L(s)$ to be a Butterworth filter, it inherently introduces complex poles to the system. Hence, from (\ref{aucorr_form}), the filtered delay cost function will exhibit a slightly underdamped response. It can then be seen that the locations of the extrema of the cost function are associated with the underdamped response of the Butterworth filter. Hence, the augmentation of $L(s)$ with $G_0(s)$ when $L(s)$ is implemented as a Butterworth filter provides a cost function similar in nature to Fig. \ref{Simple_localMinima}.

Next, to establish the sufficient condition, we need to derive a cost function $J_0$ that shares the same global minimum with the filtered cost functions where the minimum path corresponds to a decrease in $J_0$. 

To analyse this condition, we take into consideration the model error that is generated in the iterations when estimating the rational part of the system due to an inaccurate estimate of the delay. Consider we have an estimate for the delay, $\hat{\tau}$, then when estimating the rational component of (\ref{md_setting}), we have the following relationship in the Laplace domain,
\begin{equation} \label{model_err}
\hat{G}(s) = G_0(s)e^{-(\tau_0-\hat{\tau})s},
\end{equation}
where $\hat{G}(s)$ is the estimate for the rational part of the system. It is obvious that when $\hat{\tau}\neq \tau_0$, $\hat{G}(s)\neq G_0(s)$. Hence there always exists model error in estimating the rational part of the system if $\hat{\tau} \neq \tau_0$.

\begin{figure}[h]
\begin{center}
\vspace{0.3cm}
\includegraphics[width=75mm,height=52mm]{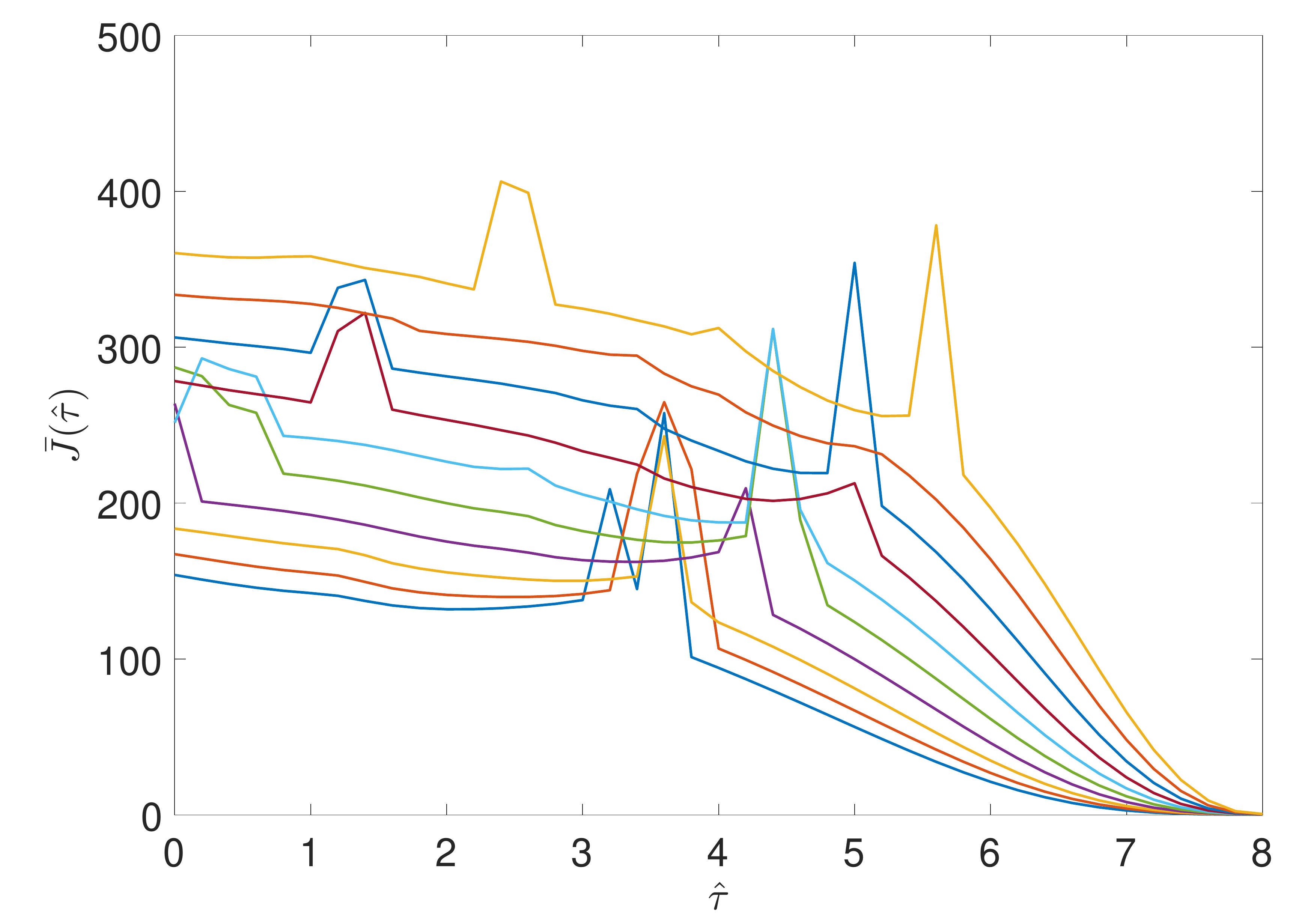}
\caption{Plot of $\bar{J}(\hat{\tau})$ for set of filters $L_k(p)$ defined as in Theorem 3.} 
\vspace{0.3cm}
\label{Decrease_setFilter} 
\end{center}
\end{figure}

We now quantify the model error to understand the impact of it on the cost function. Here, a model reduction technique based on an identification method \cite{Gu2011} is used to find a rational transfer function $G_{r}(s)$. The idea is to generate a noise-free data set from the system $G_0(s)e^{-(\tau_0-\hat{\tau})s}$, then use the SRIVC method to obtain a rational system $\hat{G}_r(s)$ that has the same model order as the true system and satisfies (\ref{model_err}). For example, Fig. \ref{Decrease_setFilter} shows the graphs of ten filtered delay cost functions after including the model error for the system $G(s) = 2e^{-8s}/(0.25s^2+0.7s+1)$. 
%It can be seen that a decrease occurs in all the filtered cost functions.

We now define the cost function $J_0$ as,
\begin{equation} \label{eq_J0}
J_0(\delta\tau) = \dfrac{1}{n_f}\sum\limits_{k=1}^{n_f} \bar{\bar{J}}_k(\delta\tau),
\end{equation}
where $\bar{\bar{J}}_k(\delta\tau) = 100\dfrac{\parallel y(t)-\hat{\bar{x}}_k(t) \parallel_2}{\parallel y(t) - \mathbb{E}\{y(t)\} \parallel_2}$, $\hat{\bar{x}}_k(t)$ is the estimated system output when the filter $L_k(s)$ is used. Note that $\bar{\bar{J}}_k(\delta\tau)$ is the normalized version of the filtered cost function, $\bar{J}_k(\delta\tau)$.

\begin{figure}[h]
\begin{center}
\vspace{0.5cm}
\includegraphics[width=75mm,height=52mm]{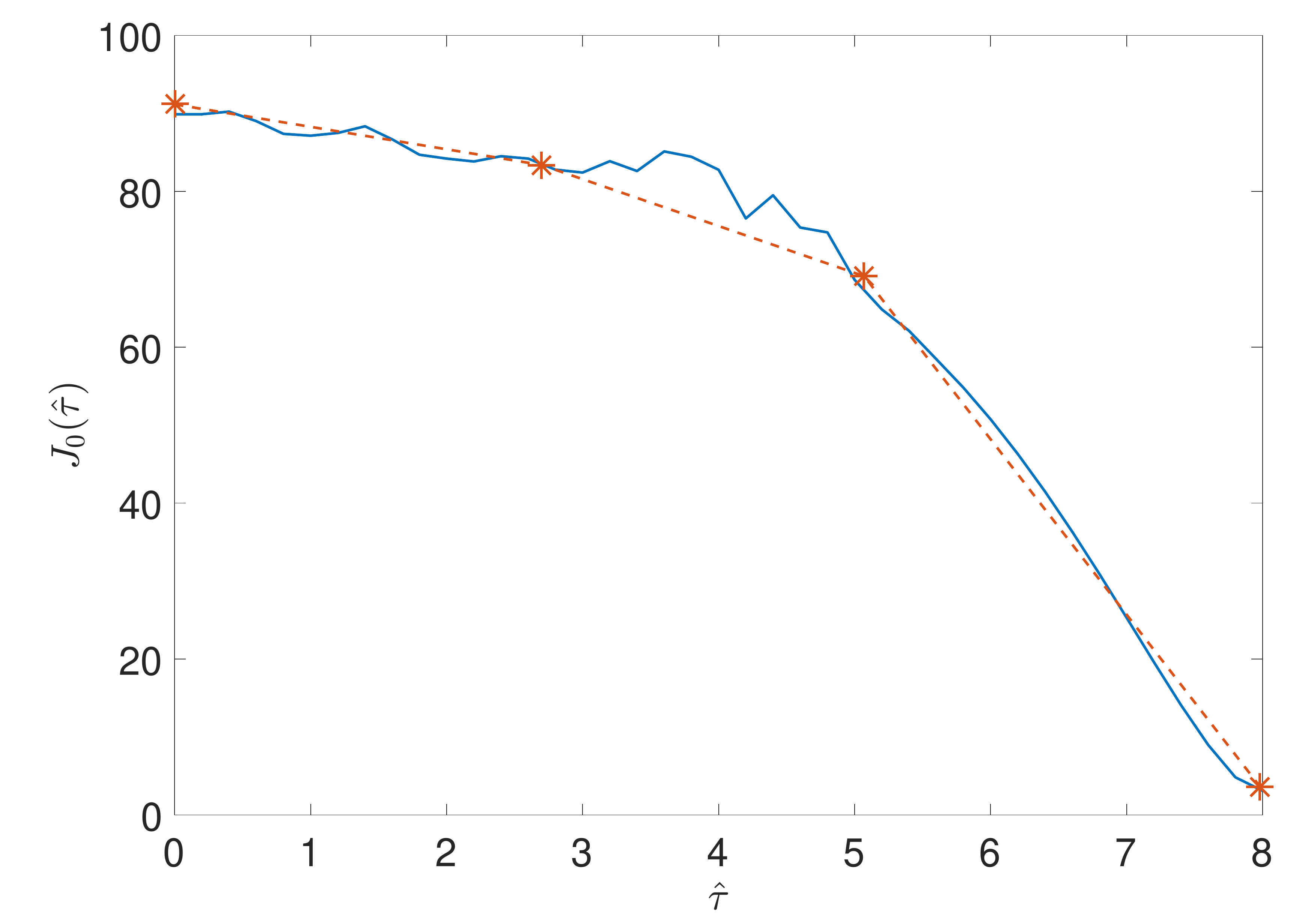}
\caption{Plot of the average normalized mean square error of the filtered cost delay function ($\ast$: switching point, solid line: $J_0$).} 
\vspace{0.3cm}
\label{Decrease_J0} 
\end{center}
\end{figure}

The averaging is used to obtain a smooth decreasing cost function $J_0$ as there are always errors in the estimation due to noise and model error. Fig. \ref{Decrease_J0} shows the average normalized mean square error for the ten filtered delay cost functions shown in Fig. \ref{Decrease_setFilter}. It can be seen that $J_0$ has the same global minimum as the ten filtered cost functions. Also shown in Fig. \ref{Decrease_J0} is the minimum path through the filtered cost function that corresponds to a decrease in $J_0$. Note that the `$\ast$' in Fig. \ref{Decrease_J0} corresponds to a switch between cost functions $J_k(\delta\tau)$.

Therefore, in summary, the \textit{M-safely redundant} optimization problem for the time delay estimation is defined as,

\begin{equation}
\underset{\tau}{\text{min}} \ \ J_0(\tau),
\end{equation}
with $J_0(\tau)$ defined as in (\ref{eq_J0}). The filter set consists of $n_f$ Butterworth filters with cut-off frequencies, chosen such that the corresponding periods are linearly spaced, that span from $1/\beta$ to 1 of the system bandwidth where $\beta \geq \dfrac{5/4-1/(5\pi^2)}{3/4-1/(3\pi^2)}$ and $n_f$ is chosen using the methodology discussed in Section 4.3. 

\subsection{Algorithm}
The proposed algorithm utilizing useful redundancy with multiple filters is described in Algorithm IV. \\
{\begin{center}
\textsc{ALGORITHM IV} 
\end{center}}
\vspace{0.3cm}
\begin{mdframed}
\textbf{Step 1. Initialization} 
\begin{enumerate}
\item
Select a set of low-pass filters as suggested in Section 4.3.
\item
Set the boundaries $\tau_{min}, \tau_{max}, \Delta \tau_{min}, \Delta \tau_{max}$\footnote{Note that $\Delta \tau_{max}$ is used to constrain the increment $\Delta\hat{\tau}$ in the time delay estimation, i.e. when $\Delta\hat{\tau}^j > \Delta \tau_{max}$, set $\Delta\hat{\tau}^j = \Delta \tau_{max}$. This is to ensure the Wolfe conditions are satisfied at each iteration in the Gauss-Newton technique \cite{Wolfe1969}.}, and the SVF cut-off frequency $\omega_c^{SVF}$. 
\item
Choose an initial value $\hat{\tau}_0$ of the time delay.
\end{enumerate}
\textbf{Step 2. Iterative estimation}

\textit{for} \textit{i=1:converge}
\begin{enumerate}
\item
\textit{for k=1:$n_{f}$}
\begin{enumerate}
\item
Choose the low-pass filter $L_k(p)$.
\item
Set the initial delay $\hat{\tau}^0_{i,k}$ as $\hat{\tau}_{i-1}$
\item
Use the SRIVC algorithm and Gauss-Newton method as described in Algorithm II (without Step 3) to estimate the time delay $\hat{\tau}_{i,k}$. 
\item
Compute $J_0(\hat{\tau}_{i,k})$ (follows (\ref{eq_J0})).
\end{enumerate}
\textit{end}
\item
Choose $\hat{\tau}_i = \underset{\hat{\tau}_{i,k}}{\text{argmin}} \ \ J_0(\hat{\tau}_{i,k})$.
\item
If $|(\hat{\tau}_i-\hat{\tau}_{i-1})/\hat{\tau}_i | \geq \epsilon$, go to Step 1, else break. \footnote{$\epsilon$ is a small value selected to obtain the desired accuracy.}
\end{enumerate}

\textit{end}
\\
\textbf{Step 3. Refine parameter estimation}

Repeat Step 2 with the low-pass filter $L(p) = 1$ .
\end{mdframed}
\vspace{0.3cm}

%%%%%%%%%%%%%%%%%%%%%%%%%%%%%%%%%%%%%%%%%%%%%%%%%%%%%%%%%%%%%%%%%%%%%%%%%%
\section{Numerical Examples}
%%%%%%%%%%%%%%%%%%%%%%%%%%%%%%%%%%%%%%%%%%%%%%%%%%%%%%%%%%%%%%%%%%%%%%%%%%

In this section, we demonstrate the effectiveness of the algorithm through numerical examples. These examples are commonly used in the literature \cite{Chen2015} \cite{Garnier2008} \cite{Rao2002} and act as defacto benchmarks to compare the performance between different CT system identification algorithms. 
\vspace{-0.1cm}
\subsection{Case 1}
\vspace{-0.1cm}
Case 1 considers the second order system used in \cite{Chen2015}, i.e.
\begin{equation} \label{sys_test_sp}
G(s) = \dfrac{2e^{-8s}}{0.25s^2+0.7s+1}.
\end{equation}
In this example, we consider two sampling schemes for the input and output data. 
\begin{enumerate}
\item
Regular (uniform) sampling: The excitation input signal is a PRBS (pseudo-random binary sequence) of maximum length. The sampling time is 50ms, the number of stages in the shift register is 10 and the number of samples, N = 4000.
\item
Irregular (nonuniform) sampling: The input excitation signal is a PRBS of maximum length. The number of stages in the shift register is 10 and the clock period is 0.5s. The input and output data are sampled at an irregular time instant $t_k$, with a sampling interval $h_k$ uniformly distributed as, $h_k \sim U[0.01, 0.09] (s)$. The realization of $h_k$ is randomized for each run. The number of samples, $N=4000$.
\end{enumerate}

The additive output disturbance is Gaussian distributed white noise with zero mean designed to give Signal to Noise Ratios (SNR) of 5dB and 15dB. The SNR is defined as,
\begin{equation} \nonumber
SNR = 10\text{log}\dfrac{P_x}{P_e},
\end{equation}
with $P_x$ and $P_e$ the average power of the noise-free output $x(t_k)$ and the disturbance noise $e(t_k)$ respectively.

The system order is assumed known a priori (order 2). The frequency $\omega_c^{SVF}$ is chosen as 1 (rad/sec) which is same as the value used in \cite{Chen2015}. To evaluate the algorithm, the time delay is initialized from different values $\hat{\tau}^0=0,1,3,5,7,9$ (similar to the experiment in \cite{Chen2015}) and also from a random value drawn from a uniform distribution $U[0, 9] (s)$.

For the two sampling schemes, 100 different data sets are generated for each noise level. The system parameters and time delay are estimated using both the algorithm described in \cite{Chen2015} and the algorithm proposed in this paper:
\begin{itemize}
\item
First we consider the algorithm developed in \cite{Chen2015}: The cut off frequency of the filter is chosen as $1/10$ of the system bandwidth as suggested in \cite{Chen2015}. For regular data, the command \textit{`tdsrivc'} from the CONTSID toolbox \cite{Padilla15} is used to estimate the time delay and the system parameters. 
\item
For the multiple filtering algorithm utilizing the useful redundancy technique as proposed in this paper, the number of filters is set to 10 and the cut off frequencies, chosen based on linearly spaced periods, of the filters span from $1/10$ to $1$ of the system bandwidth. The order of the Butterworth filter is set to 10.
\end{itemize}

For both algorithms, the time delay boundaries $\tau_{min}, \tau_{max}$ are set to $0$ and $15$ seconds respectively. The maximum number of iterations for the Gauss-Newton method (Step 2 in Algorithm II and IV) is set to 10. The threshold, $\epsilon$, to determine convergence is set to $10^{-3}$.

%\begin{table*}[t]
%\centering
%\caption{Global convergence as a percentage of 100 Monte Carlo simulations for case 1}
%\begin{tabular}{|p{2.0cm}|p{0.5cm}|p{0.5cm}|p{0.5cm}|p{0.5cm}|p{0.5cm}|p{0.5cm}|p{0.5cm}|p{0.5cm}|p{0.5cm}|p{0.5cm}|p{0.5cm}|p{0.5cm}|p{0.45cm}|p{0.45cm}|p{0.55cm}|}
%\hline
%\multirow{4}{*}{\parbox{2.0cm}{\centering Sampling scheme}} & \multirow{4}{*}{SNR} & \multicolumn{6}{l|}{\centering Existing method \cite{Chen2015}} & \multicolumn{6}{l|} {\centering Proposed method} \\ [2ex] \cline{3-14} 
%\multicolumn{1}{|l|}{} & & 0s & 1s & 3s & 5s & 7s & 9s & 0s & 1s & 3s & 5s & 7s & 9s \\[0.5ex] \hline
%\multirow{2}{*}{\parbox{2.0cm}{\centering Regular sampled data}} & 5dB & 65\% & 100\% & 100\% & 97\% & 100\% & 100\% & 100\% & 100\% & 100\% & 100\% & 100\% & 100\% \\[0.5ex] \cline{2-14} 
%\multicolumn{1}{|l|}{} & 15dB & 69\% & 89\% & 92\% & 99\% & 100\% & 100\% & 100\% & 100\% & 100\% & 100\% & 100\% & 100\% \\[0.5ex] \hline
%\multirow{2}{*}{\parbox{2.0cm}{\centering Irregular sampled data}} & 5dB & 1\% & 1\% & 0\% & 17\% & 100\% & 90\% & 100\% & 100\% & 100\% & 100\% & 100\% & 100\% \\[0.5ex] \cline{2-14} 
%\multicolumn{1}{|l|}{} & 15dB & 1\% & 0\% & 7\% & 18\% & 100\% & 89\% & 100\% & 100\% & 100\% & 100\% & 100\% & 100\% \\[0.5ex] \hline
%\end{tabular}
%\label{table_gcTD1}
%\end{table*}

\begin{table*}[t]
\centering
\caption{Global convergence as a percentage of 100 Monte Carlo simulations for case 1}
\begin{tabular}{|p{2.0cm}|p{1.8cm}|p{0.53cm}|p{0.53cm}|p{0.53cm}|p{0.53cm}|p{0.53cm}|p{0.53cm}|p{0.53cm}|p{0.53cm}|p{0.53cm}|p{0.53cm}|p{0.53cm}|p{0.52cm}|p{0.55cm}|p{0.8cm}|}
\hline
\multirow{2}{*}{\parbox{2.0cm}{\centering Sampling scheme}} & \multicolumn{1}{l|}{} & \multicolumn{6}{l|}{\centering Existing method \cite{Chen2015}} & \multicolumn{6}{l|} {\centering Proposed method} \\ [0.4ex] \cline{2-14} 
\multicolumn{1}{|l|}{} & Initial delay & 0s & 1s & 3s & 5s & 7s & 9s & 0s & 1s & 3s & 5s & 7s & 9s \\[0.4ex] \hline
\multirow{2}{*}{\parbox{2.0cm}{\centering Regular sampled data}} & SNR = 5dB & 65\% & 100\% & 100\% & 97\% & 100\% & 100\% & 100\% & 100\% & 100\% & 100\% & 100\% & 100\% \\[0.4ex] \cline{2-14} 
\multicolumn{1}{|l|}{} & SNR = 15dB & 69\% & 89\% & 92\% & 99\% & 100\% & 100\% & 100\% & 100\% & 100\% & 100\% & 100\% & 100\% \\[0.4ex] \hline
\multirow{2}{*}{\parbox{2.0cm}{\centering Irregular sampled data}} & SNR = 5dB & 1\% & 1\% & 0\% & 17\% & 100\% & 90\% & 100\% & 100\% & 100\% & 100\% & 100\% & 100\% \\[0.4ex] \cline{2-14} 
\multicolumn{1}{|l|}{} & SNR = 15dB & 1\% & 0\% & 7\% & 18\% & 100\% & 89\% & 100\% & 100\% & 100\% & 100\% & 100\% & 100\% \\[0.4ex] \hline
\end{tabular}
\label{table_gcTD1}
\end{table*}

We consider the estimated time delay to be the global minimum when the relative error, $\epsilon_r$, is less than $1\%$, i.e. 
\begin{equation} \label{err_rel}
\epsilon_r = \dfrac{| \hat{\tau}-\tau_0 |}{\tau_0}\times 100\%,
\end{equation}
where $\hat{\tau}$ and $\tau_0$ are the estimated delay and the true system delay respectively.

\begin{figure} [h]
\begin{center}
\includegraphics[width=75mm,height=52mm]{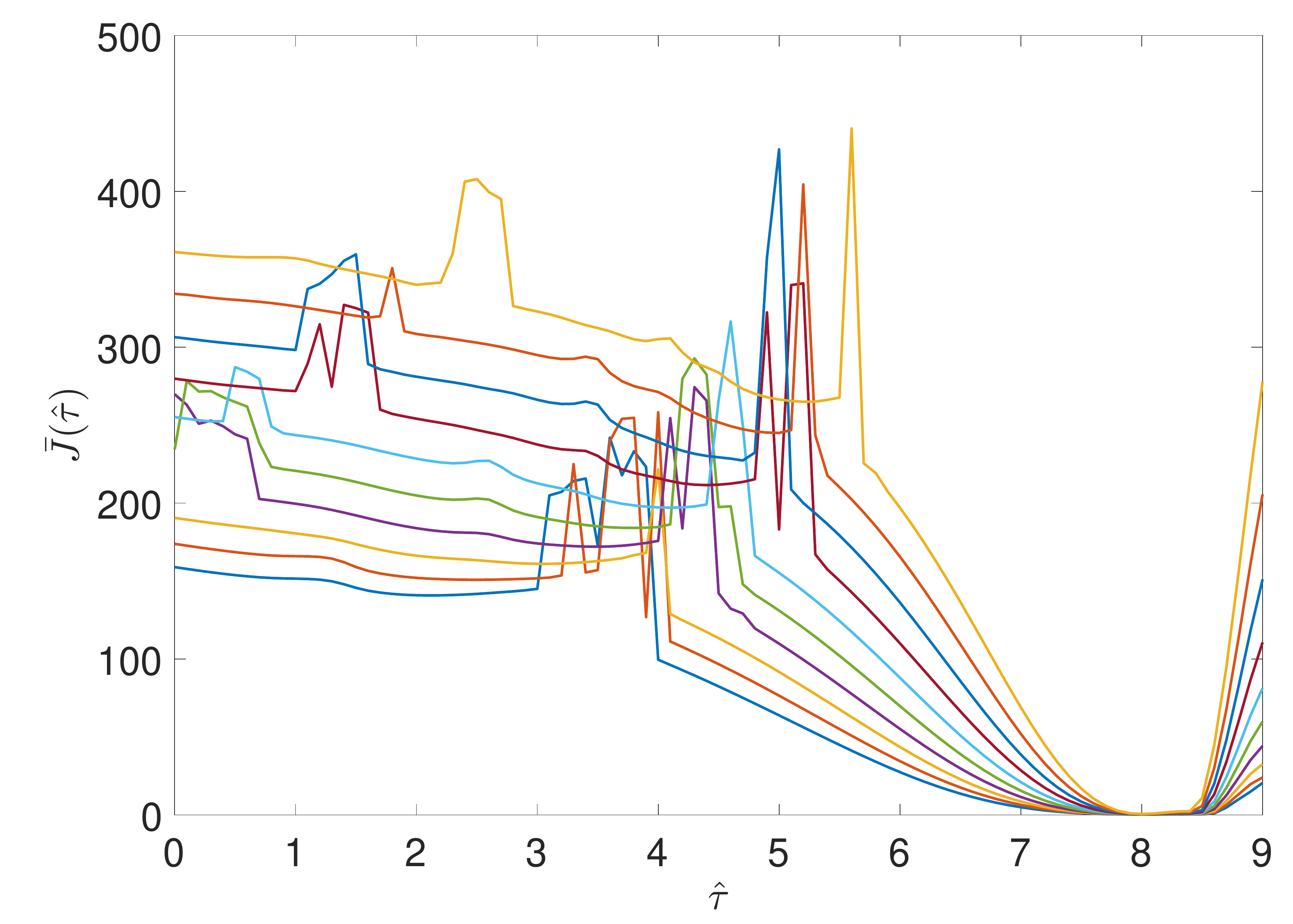}
\caption{Case 1: Cost functions of the 10 filters.} %
\label{CostFunc_sp} 
\end{center}
\end{figure}

Fig. \ref{CostFunc_sp} contains plots of the 10 cost functions corresponding to the filters used in the useful redundancy method. It can be seen that they all share the same global minimum, however they all possess different local minima. 

\begin{figure} [h]
\begin{center}
\includegraphics[width=75mm,height=52mm]{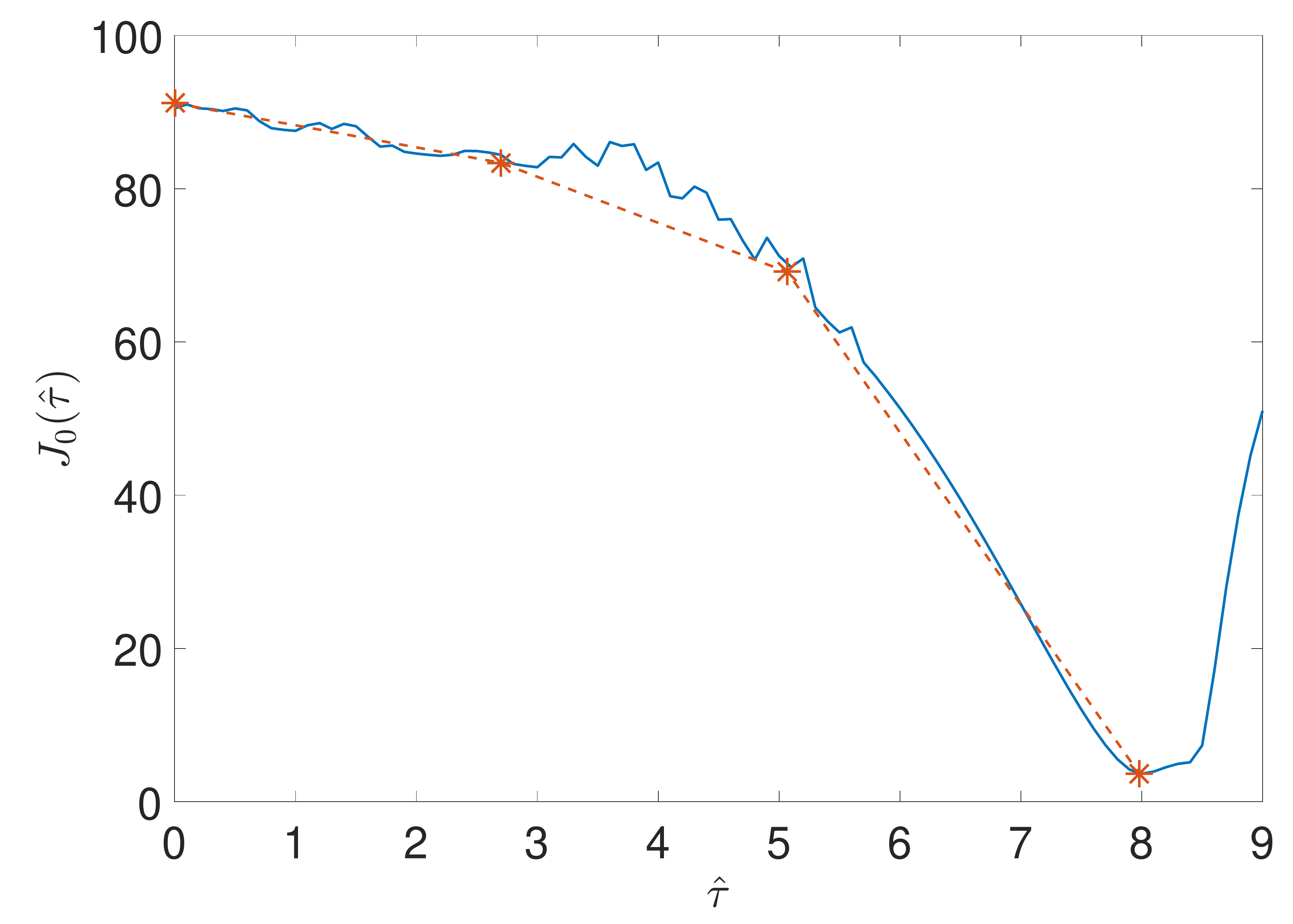}
\caption{Case 1: The original cost function and time delay estimate trajectory for a single dataset ($\ast$: switching point, solid line: $J_0$).} 
\label{SolverPath1} 
\end{center}
\end{figure}

Fig. \ref{SolverPath1} presents a plot of $J_0$, for a single dataset, showing the trajectory of the time delay estimate when the sampled data is regular and the SNR is 15dB. The algorithm is started with the initial value of the time delay, $\hat{\tau}^0 = 0$. The `*' in Fig. \ref{SolverPath1} represents the algorithm switching to another solver path. The algorithm stops when the condition in Algorithm IV is satisfied for the threshold $\epsilon = 10^{-3}$, i.e. when the relative error in the time delay between consecutive iterations is $10^{-3}$. 

Next we demonstrate the effectiveness of the proposed algorithm by showing the global convergence of 100 datasets with different initial values for the delay. The results are provided in Tables \ref{table_gcTD1} and \ref{table_random1} where they are also compared to results from the existing method \cite{Chen2015}. Note that the results presented here for the existing method are different to those achieved in \cite{Chen2015} as we chose a larger bound on the delay. In \cite{Chen2015}, $\tau_{max}$ was set to 10 seconds while here, it is set to 15 seconds.

\begin{table}[h]
\centering
\caption{Global convergence as a percentage of 100 Monte Carlo simulations for case 1 with a randomly selected initial delay}
\begin{tabular}{|p{2cm}|p{0.5cm}|p{1.8cm}|p{1.3cm}|p{0.5cm}|p{0.5cm}|p{0.5cm}|p{0.5cm}|p{0.5cm}|p{0.5cm}|p{0.5cm}|p{0.5cm}|p{0.5cm}|p{0.45cm}|p{0.45cm}|p{0.45cm}|}
\hline
\multirow{2}{*}{\parbox{1.2cm}{Sampling scheme}} & \multirow{2}{*}{SNR} & \multirow{2}{*}{\parbox{1.8cm}{\centering Existing method \cite{Chen2015}}} & \multirow{2}{*}{\parbox{1.3cm}{\centering Proposed method}} \\[2.4ex] \hline
\multirow{2}{*}{\parbox{1.9cm}{\centering Regular sampled data}}  & 5dB & \centering 90\% & 100\% \\[0.4ex] \cline{2-4} 
\multicolumn{1}{|l|}{} & 15dB & \centering 97\% & 100\% \\[0.4ex] \hline
\multirow{2}{*}{\parbox{1.9cm}{\centering Irregular sampled data}} & 5dB & \centering 40\% & 100\% \\[0.4ex] \cline{2-4} 
\multicolumn{1}{|l|}{} & 15dB & \centering 41\% & 100\% \\[0.4ex] \hline
\end{tabular}
\label{table_random1}
\end{table}

The results presented in Tables \ref{table_gcTD1} and \ref{table_random1} show that the proposed method, which utilizes the useful redundancy technique, performs much better than the existing method \cite{Chen2015} irrespective of the initial delay. When the initial delays are close to the true system delay, the global convergence percentage is high, i.e. approximately 100\% for both algorithms. However, when the initial delays are poorly selected, the global convergence percentage is low for the existing method, e.g. with the initial delay $= 0$, the percentage convergence for the existing method \cite{Chen2015} is less than 1\% while the proposed method achieves 100\% global percentage convergence for both SNRs and the two sampling schemes.

\subsection{Case 2: Rao-Garnier test system}

\begin{table*}[t]
\centering
\caption{Global convergence as a percentage of 100 Monte Carlo simulations for Case 2}
\begin{tabular}{|p{2.0cm}|p{1.8cm}|p{0.53cm}|p{0.53cm}|p{0.53cm}|p{0.53cm}|p{0.53cm}|p{0.53cm}|p{0.53cm}|p{0.53cm}|p{0.53cm}|p{0.53cm}|p{0.53cm}|p{0.53cm}|p{0.5cm}|p{0.5cm}|}
\hline
\multirow{2}{*}{\parbox{2.0cm}{\centering Sampling scheme}} & \multicolumn{1}{l|}{} & \multicolumn{6}{l|}{\centering Existing method \cite{Chen2015}} & \multicolumn{6}{l|} {Proposed method} \\ [0.4ex] \cline{2-14} 
\multicolumn{1}{|l|}{} & Initial delay & 0s & 1s & 3s & 5s & 7s & 9s & 0s & 1s & 3s & 5s & 7s & 9s \\[0.4ex] \hline
\multirow{2}{*}{\parbox{2.0cm}{\centering Regular sampled data}} & SNR = 5dB & 1\% & 52\% & 22\% & 43\% & 76\% & 68\% & 99\% & 99\% & 99\% & 99\% & 99\% & 99\% \\[0.4ex] \cline{2-14} 
\multicolumn{1}{|l|}{} & SNR = 15dB & 0\% & 60\% & 6\% & 28\% & 70\% & 69\% & 100\% & 100\% & 100\% & 100\% &  100\% & 100\% \\[0.4ex] \hline
\multirow{2}{*}{\parbox{2.0cm}{\centering Irregular sampled data}} & SNR = 5dB & 8\% & 16\% & 23\% & 26\% & 66\% & 48\% & 100\% & 100\% & 100\% & 100\% & 100\% & 100\% \\[0.4ex] \cline{2-14} 
\multicolumn{1}{|l|}{} & SNR = 15dB & 9\% & 10\% & 19\% & 27\% & 55\% & 50\% & 100\% & 100\% & 100\% & 100\% & 100\% & 100\% \\[0.4ex] \hline
\end{tabular}
\label{table_gcTD}
\end{table*}

In this section, we consider a system based on the Rao-Garnier continuous time benchmark \cite{Garnier2008}, 
\begin{equation} \nonumber
G(s) = \dfrac{(-6400s+1600)e^{-8s}}{s^4+5s^3+408s^2+416s+1600}.
\end{equation}

This system is linear, non-minimum phase with complex poles and a time delay. As with case 1, the experiment is conducted using both regular and irregular sampling schemes. For each sampling scheme, the input excitation signal is a PRBS of maximum-length. In the regular sampling case, the sample time is 10ms, the number of stages in the shift register is 10, the number of samples, N = 8000. For the irregular sampling, the number of stages in the shift register is 10 and the clock period is 0.5s. The input and output data are sampled at an irregular time instant $t_k$, where the sampling interval $h_k$ is uniformly distributed as, $h_k \sim U[0.01, 0.05] (s)$. The number of samples, $N=4500$.

We consider two SNR levels: 5dB and 15dB. The frequency $\omega_c^{SVF}$ is chosen as 25 rad/sec, which is approximately the bandwidth of the system. Again, the initial values for the time delay, $\hat{\tau}^0$, are selected as $0, 1, 3, 5, 7, 9$ as well as a random value from the uniform distribution $U[0, 9] (s)$. The delay boundaries $\tau_{min}, \tau_{max}$ are set to $0$ and $15$ sec respectively. 

For the existing method, the cut off frequency of the filter is chosen as $1/10$ of the system bandwidth. For the proposed method, the cut off frequencies of the filters, chosen based on linearly spaced periods, span from $1/10$ to $1$ of the system bandwidth and the number of filters used is 15. The order of the Butterworth filters is set to 10. The maximum number of iterations for the Gauss-Newton method (Step 2 in Algorithm II and IV) is 10. The threshold, $\epsilon$, for convergence is $10^{-3}$.

For comparison, 100 different data sets are generated for each noise level of each sampling scheme.

\begin{figure} [h]
\begin{center}
\includegraphics[width=75mm,height=52mm]{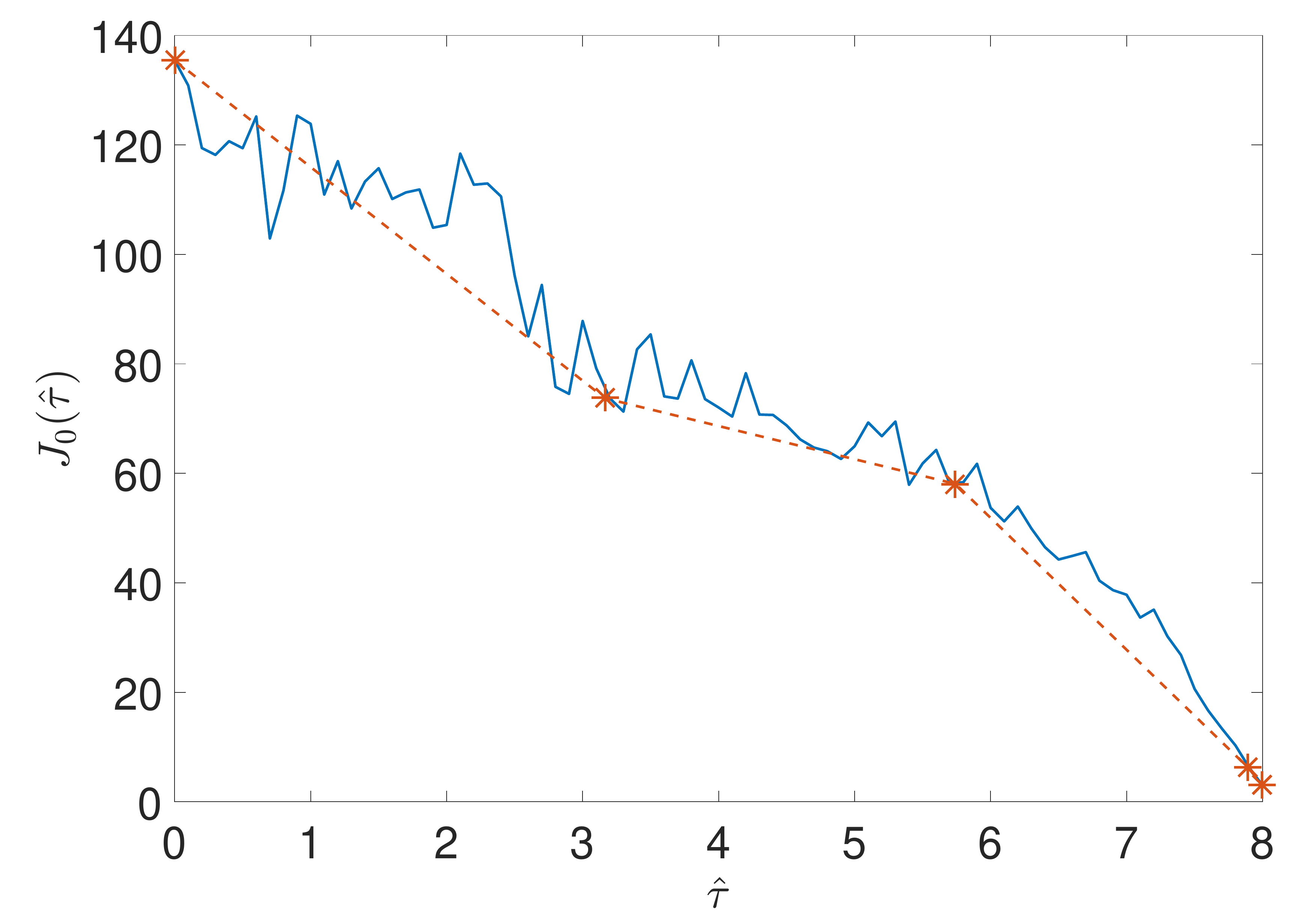}
\caption
{Case 2: The original cost function and time delay estimate trajectory for a single data set ($\ast$: switching point, solid line: $J_0$)}
\label{SolverPath} 
\end{center}
\end{figure}

A graph showing the trajectory of a time delay estimate is plotted in Fig. \ref{SolverPath} for a single data set. Switching between cost functions is indicated by '$*$'. Note that, if only one filter were to be used, there is a high chance of convergence to a local minima. However, by using the useful redundancy technique, the proposed algorithm can avoid the local minima and converge to the global minimum.

\begin{table}[h]
\centering
\caption{Global convergence as a percentage of 100 Monte Carlo simulations for case 2 with a randomly selected initial delay}
\begin{tabular}{|p{1.8cm}|p{0.5cm}|p{1.8cm}|p{1.3cm}|p{0.5cm}|p{0.5cm}|p{0.5cm}|p{0.5cm}|p{0.5cm}|p{0.5cm}|p{0.5cm}|p{0.5cm}|p{0.5cm}|p{0.45cm}|p{0.45cm}|p{0.45cm}|}
\hline
\multirow{2}{*}{\parbox{1.6cm}{\centering Sampling scheme}} & \multirow{2}{*}{SNR} & \multirow{2}{*}{\parbox{2cm}{\centering Existing method \cite{Chen2015}}} & \multirow{2}{*}{\parbox{1.3cm}{\centering Proposed method}} \\[2ex] \hline
\multirow{2}{*}{\parbox{1.8cm}{\centering Regular sampled data}}  & 5dB & \centering 52\% & 99\% \\[0.4ex] \cline{2-4} 
\multicolumn{1}{|l|}{} & 15dB & \centering 46\% & 100\% \\[0.4ex] \hline
\multirow{2}{*}{\parbox{1.8cm}{\centering Irregular sampled data}} & 5dB & \centering 48\% & 100\% \\[0.4ex] \cline{2-4} 
\multicolumn{1}{|l|}{} & 15dB & \centering 21\% & 100\% \\[0.4ex] \hline
\end{tabular}
\label{table_random2}
\end{table}

The numerical results of the experiments are provided in Table \ref{table_gcTD} and \ref{table_random2}. We can see that for both sampling schemes, the proposed method based on the useful redundancy technique performs much better as compared to the existing method. For all the initial values of delay used in this experiment, the existing method \cite{Chen2015} never convergences 100\% to the global minimum. However the proposed method utilizing useful redundancy still achieves a very high global convergence percentage, i.e. mostly 100\% with any initial delay for both SNRs and sampling schemes. 

%%%%%%%%%%%%%%%%%%%%%%%%%%%%%%%%%%%%%%%%%%%%%%%%%%%%%%%%%%%%%%%%%%%%%%%%%%
\section{Conclusion}
%%%%%%%%%%%%%%%%%%%%%%%%%%%%%%%%%%%%%%%%%%%%%%%%%%%%%%%%%%%%%%%%%%%%%%%%%%

The paper presents a new algorithm to estimate the parameters and time delay of a continuous-time system from regularly and irregularly sampled data. The idea is based on Instrumental Variable methods and employing the useful redundancy technique to enhance the global convergence by generating multiple cost functions by filtering the data several times. The paper also develops theoretical results related to the minima locations of the filtered delay cost function and the choice of filters to ensure the algorithm can converge to a global minimum. Numerical results show a significant improvement in the global convergence rate of the time delay estimation as compared to existing methods irrespective of the SNR.

%%%%%%%%%%%%%%%%%%%%%%%%%%%%%%%%%%%%%%%%%%%%%%%%%%%%%%%%%%%%%%%%%%%%%%%%%%%%%%%%

\bibliographystyle{plain}        % Include this if you use bibtex 
\bibliography{library}           % and a bib file to produce the 

\appendix
\section{Appendix}

\vspace{-0.2cm}
\subsection{Proof of Theorem 1}
\begin{proof}
From (\ref{J_filter2}), when $\lambda=0$, $\bar{J}(\delta\tau)$ becomes,
\begin{equation}
\begin{aligned}
& \bar{J}(\delta\tau) = \dfrac{1}{\pi}\int_{-\infty}^{\infty} \! (1-\cos(\omega\delta\tau)) \vert G_0(j\omega) L(j\omega)\big\vert^2 \, \mathrm{d}\omega.
\end{aligned}
\end{equation}
$\forall \delta\tau \in \mathbb{R},\ \omega \geq 0$, we have
\begin{equation} \nonumber
\begin{aligned}
& 1-\cos(\omega\delta\tau) \geq 0, \ \ \forall \delta\tau, \\
& \vert G_0(j\omega) L(j\omega)\big\vert^2 \geq 0, \ \ \forall \omega.
\end{aligned}
\end{equation}
Therefore, $\bar{J}(\delta\tau) \geq 0,\ \forall \delta\tau \in \mathbb{R},\ \omega \geq 0$. Note that $\bar{J}(0) = 0$, hence we have $\bar{J}(\delta\tau) \geq \bar{J}(0),\ \forall \delta\tau \in \mathbb{R},\ \omega \geq 0$. The equality occurs when $L(j\omega) = 0,\ \forall \omega$ or $G_0(j\omega) = 0,\ \forall \omega$ or $\cos(\omega\delta\tau)) = 1,\ \forall \omega$. As the system $G_0(p)$ and the filter $L(p)$ are not 0 (condition of Theorem 1), the equality only occurs when $\cos(\omega\delta\tau)) = 1,\ \forall \omega$, hence $\delta\tau = 0$. \hfill $\square$
\end{proof}

\vspace{-0.2cm}
\subsection{Proof of Theorem 2}
\begin{proof}
If $L(s)$ is selected such that $L(s)G_0(s)$ is an ideal low-pass filter with cut-off frequency $\omega_c$, then $\bar{J}(\delta\tau)$ becomes,
\begin{equation} \label{J_eq}
\begin{aligned}
\bar{J}(\delta\tau) & = \dfrac{1}{2\pi}\int_{-\omega_c}^{\omega_c} \! \vert 1-e^{-j\omega\delta\tau} \vert^2 \mathrm{d}\omega + C\\
& = \dfrac{2\omega_c}{\pi}-\dfrac{2\sin(\omega_c\delta\tau)}{\pi\delta\tau} + C,
\end{aligned}
\end{equation}
where $C = \dfrac{1}{2\pi}\int_{-\infty}^{\infty} \! \lambda \big\vert L(j\omega)\big\vert^2 \mathrm{d}\omega$.

The minima and maxima of $\bar{J}(\delta\tau)$ occur at the roots of,
\begin{equation} \label{J_delta}
\begin{aligned} = \dfrac{\mathrm{d}\bar{J}(\delta\tau)}{\mathrm{d}\delta\tau} = -\dfrac{2}{\pi}\dfrac{\omega_c\cos(\omega_c\delta\tau)\delta\tau-\sin(\omega_c\delta\tau)}{\delta\tau^2}.
\end{aligned}
\end{equation}
From (\ref{J_delta}), we can see that the extrema of $\bar{J}(\delta\tau)$ are also the extrema of the function $sinc(\omega_c\delta\tau)$. However, the minima of $\bar{J}(\delta\tau)$ will be the maxima of $sinc(\omega_c\delta\tau)$ and the maxima of $\bar{J}(\delta\tau)$ will be the minima of $sinc(\omega_c\delta\tau)$. Note that for the function $sinc(\omega_c\delta\tau)$, the locations of the $i^{th}$ ($i \geq 1$) positive extremum $\tilde{\delta}\tau$ can be approximated by,
\begin{equation} \label{extremum_sinc}
\begin{aligned}
\tilde{\delta}\tau_i & \simeq (i+\dfrac{1}{2})\dfrac{\pi}{\omega_c} - \dfrac{1}{(i+\dfrac{1}{2})\pi \omega_c} \\
& \simeq \dfrac{2i+1}{4}T_c - \dfrac{T_c}{(2i+1)\pi^2}.
\end{aligned}
\end{equation}
where $T_c$ is the corresponding period of the cut-off frequency $\omega_c$. The minima occurs when $i$ is odd and the maxima occurs when $i$ is even. 

Therefore, for $\bar{J}(\delta\tau)$, the positive extrema can also be approximated using (\ref{extremum_sinc}), the only difference is that the minima of $\bar{J}(\delta\tau)$ occurs when $i$ is even and the maxima occurs when $i$ is odd. \hfill $\square$
\end{proof}

\vspace{-0.2cm}
\subsection{Proof of Theorem 3}
\begin{proof}
Here we use the positive extrema locations of the time delay cost function in (\ref{eq_minloc}). Since the filtered time delay cost function $\bar{J}(\delta\tau)$ is an even function, we only need to prove the theorem for $\delta\tau_0 > 0$.

Denote $\tilde{\delta}\tau_{i,k}^{min}$ as the $i^{th}$ positive minimum of the filtered time delay cost function $\bar{J}(\delta\tau)$ generated by the filter $L_k(p)$; $\tilde{\delta}\tau_{i,k}^{max}$ as the $i^{th}$ positive maximum of the filtered time delay cost function $\bar{J}(\delta\tau)$ generated by the filter $L_k(p)$. From Theorem \ref{theorem_minloc}, they can be computed by,
\begin{equation} \label{eq_t3_exloc}
\begin{aligned}
\tilde{\delta}\tau_{i,k+1}^{min} & = \dfrac{4i+1}{4}T_{c,k+1} - \dfrac{T_{c,k+1}}{(4i+1)\pi^2}, \\
& = (\dfrac{4i+1}{4} - \dfrac{1}{(4i+1)\pi^2})\Big(\beta-k\dfrac{\beta-1}{n_f-1}\Big)T_{bw} \\
\tilde{\delta}\tau_{i,k+1}^{max} & = \dfrac{4i-1}{4}T_{c,k+1} - \dfrac{T_{c,k+1}}{(4i-1)\pi^2}, \\
& = (\dfrac{4i-1}{4} - \dfrac{1}{(4i-1)\pi^2})\Big(\beta-k\dfrac{\beta-1}{n_f-1}\Big)T_{bw}. \\
\end{aligned}
\end{equation}
From (\ref{eq_t3_exloc}), for any filter $L_k(p)$, we have the following property,
\begin{equation} \label{eq_t3_exlocprop}
\begin{aligned}
 \tilde{\delta}\tau_{i,k}^{max} < \tilde{\delta}\tau_{i,k}^{min} < \tilde{\delta}\tau_{i+1,k}^{max}, \ \ \forall i \in \mathbb{Z}^+.
\end{aligned}
\end{equation}

Consider the two following cases:
\vspace{0.4cm} \\
\textbf{\textit{Case 1:}} $\delta\tau_0 \leq \tilde{\delta}\tau_{1,1}^{max}$, then from (\ref{eq_t3_exlocprop}), we have, $\xi(L_1(p),\delta\tau_0) = 0$, which is smaller than $\delta\tau_0$ as $\delta\tau_0 \neq 0$. Therefore, $\xi(L_1(p),\delta\tau_0) < \delta\tau_0$. \hfill $\square$ 
\vspace{0.4cm} \\
\textbf{\textit{Case 2:}} $\delta\tau_0 > \tilde{\delta}\tau_{1,1}^{max}$. Denote $i_0$ such that,
\begin{equation}
\tilde{\delta}\tau_{i_0-1,1}^{max} < \delta\tau_0 \leq \tilde{\delta}\tau_{i_0,1}^{max}, \ \ \ \ i_0 \in \mathbb{N}, i_0 \geq 2.
\end{equation}
Now consider a filter set $L_k(p), k = \overline{1,n_f}$ defined as in Theorem \ref{theorem_nf}, with,
\begin{equation} \label{eq_t3_nf}
n_f \geq (1/M+\beta-2)/(1/M-1), \ \ \ n_f \in \mathbb{N}, 
\end{equation}
where 
\begin{equation} \label{eq_t3_M}
M = \dfrac{4i_0-3-4/((4i_0-3)\pi^2)}{4i_0-1-4/((4i_0-1)\pi^2)}.
\end{equation}
We need to prove the filter set satisfies the condition in (\ref{eq_t3_cond}). First, we prove that the filter set $L_k(p), k = \overline{1,n_f}$ has the following property,
\begin{equation} \nonumber
\tilde{\delta}\tau_{i+1,k+1}^{max} \geq \tilde{\delta}\tau_{i,k}^{min},
\end{equation}
where $i = \overline{1,i_0-1}, k = \overline{1,n_f-1}$.
\\
From (\ref{eq_t3_nf}), 
\begin{equation} \label{eq_t3_1}
\begin{aligned}
& n_f \geq \dfrac{1/M+\beta-2}{1/M-1} \\
& \ \ \ \ \Rightarrow  M \leq 1 - \dfrac{(\beta-1)/(n_f-1)}{\beta - (n_f-2)(\beta-1)/(n_f-1)}.
\end{aligned}
\end{equation}
Note that, as $\beta > 1$,
\begin{equation} \label{eq_t3_2}
\begin{aligned}
1 - \dfrac{(\beta-1)/(n_f-1)}{\beta} & > 1 - \dfrac{(\beta-1)/(n_f-1)}{\beta - (\beta-1)/(n_f-1)} \\
& > 1 - \dfrac{(\beta-1)/(n_f-1)}{\beta - (n_f-2)(\beta-1)/(n_f-1)}.
\end{aligned}
\end{equation}
From (\ref{eq_t3_1}) and (\ref{eq_t3_2}), 
\begin{equation} \label{eq_t3_2aa}
\begin{aligned}
& 1 - \dfrac{(\beta-1)/(n_f-1)}{\beta - (k-1)(\beta-1)/(n_f-1)} \geq M, 
\end{aligned}
\end{equation}
where $k = \overline{1,n_f-1}$,
hence, combining (\ref{eq_t3_M}) and (\ref{eq_t3_2aa}),
\begin{equation} \label{eq_t3_2a}
\begin{aligned}
& \dfrac{\beta-k(\beta-1)/(n_f-1)}{\beta - (k-1)(\beta-1)/(n_f-1)} \geq \dfrac{4i_0-3-4/((4i_0-3)\pi^2)}{4i_0-1-4/((4i_0-1)\pi^2)}, 
\end{aligned}
\end{equation}
where $k = \overline{1,n_f-1}.$ \\
Next, consider the function 
\begin{equation} \label{eq_fx}
f(x) = \dfrac{4x+1-4/((4x+1)\pi^2)}{4x+3-4/((4x+3)\pi^2)}.
\end{equation}
Taking the derivative w.r.t. $x$ of (\ref{eq_fx}),

\begin{equation} \nonumber
\begin{aligned}
& \dfrac{\partial f(x)}{\partial x} \\
& = \dfrac{8+\dfrac{16}{(4x+1)\pi^2}+\dfrac{16(4x+3)}{(4x+1)^2\pi^2}+\dfrac{64/\pi^4}{(4x+1)(4x+3)^2}}{\Big((4x+3)-4/((4x+3)\pi^2)\Big)^2} \\
& \ \ \ \ \ \ \ \ - \dfrac{\dfrac{16}{(4x+3)\pi^2}+\dfrac{16(4x+1)}{(4x+3)^2\pi^2}+\dfrac{64/\pi^4}{(4x+1)^2(4x+3)}}{\Big((4x+3)-4/((4x+3)\pi^2)\Big)^2} 
\end{aligned}
\end{equation}
Now $\forall x \geq 1$,
\begin{equation} \nonumber
\begin{aligned}
& \dfrac{\partial f(x)}{\partial x} \\
& > \dfrac{8 - \dfrac{16}{(4x+3)\pi^2} - \dfrac{16(4x+1)}{(4x+3)^2\pi^2} - \dfrac{64/\pi^4}{(4x+1)^2(4x+3)}}{\Big((4x+3)-4/((4x+3)\pi^2)\Big)^2} \\
& > \dfrac{8 - \dfrac{16}{(4x+3)\pi^2} - \dfrac{16}{(4x+3)\pi^2} - \dfrac{64/\pi^4}{(4x+1)^2(4x+3)}}{\Big((4x+3)-4/((4x+3)\pi^2)\Big)^2} \\
& > \dfrac{8 - 2\dfrac{16}{7\pi^2} - \dfrac{64}{175\pi^4}}{\Big((4x+3)-4/((4x+3)\pi^2)\Big)^2} \ > 0. 
\end{aligned}
\end{equation}
Hence, $f(x)$ is an increasing function for $x \geq 1$, so we have, $f(i_0-1) \geq f(i), i=\overline{1,i_0-1}$, or, 
\begin{equation} \label{eq_t3_2b}
\begin{aligned}
& \dfrac{4i_0-3-4/((4i_0-3)\pi^2)}{4i_0-1-4/((4i_0-1)\pi^2)} \\
& \ \ \ \ \ \ \geq \dfrac{4i+1-4/((4i+1)\pi^2)}{4i+3-4/((4i+3)\pi^2)}, \ \ \ \ i=\overline{1,i_0-1}.
\end{aligned}
\end{equation}
From (\ref{eq_t3_2a}) and (\ref{eq_t3_2b}),
\begin{equation} \nonumber
\begin{aligned}
& \dfrac{\beta-k(\beta-1)/(n_f-1)}{\beta - (k-1)(\beta-1)/(n_f-1)} \geq \dfrac{4i+1-4/((4i+1)\pi^2)}{4i+3-4/((4i+3)\pi^2)}, 
\end{aligned}
\end{equation}
where $i = \overline{1,i_0-1}, k = \overline{1,n_f-1}$,
which can be rewritten as, 
\begin{equation} \label{eq_t3_3}
\begin{aligned}
& \Big(\dfrac{4i+3}{4}-\dfrac{1}{(4i+3)\pi^2}\Big) \Big(\beta-k\dfrac{\beta-1}{n_f-1}\Big) \\
& \geq \Big(\dfrac{4i+1}{4}-\dfrac{1}{(4i+1)\pi^2}\Big) \Big(\beta-(k-1)\dfrac{\beta-1}{n_f-1}\Big),
\end{aligned}
\end{equation}
for $i = \overline{1,i_0-1}, k = \overline{1,n_f-1}$. \\
From (\ref{eq_t3_3}) and (\ref{eq_t3_exloc}), we conclude,
\begin{equation} \label{eq_t3_cond1}
\tilde{\delta}\tau_{i+1,k+1}^{max} \geq \tilde{\delta}\tau_{i,k}^{min}, \ \ \ \ i = \overline{1,i_0-1}, k = \overline{1,n_f-1}.
\end{equation}

Now, we will prove that $\tilde{\delta}\tau_{i,1}^{max} \geq \tilde{\delta}\tau_{i,n_f}^{min}$, $i = \overline{1,i_0}$. Similar to the previous analysis, it's trivial to show that,
\begin{equation} \nonumber
\begin{aligned}
\dfrac{5-4/(5\pi^2)}{3-4/(3\pi^2)} \geq \dfrac{4i+1-4/((4i+1)\pi^2)}{4i-1-4/((4i-1)\pi^2)},\ \forall i \in \mathbb{Z}^+ \\
\end{aligned}
\end{equation}
As $\beta \geq \dfrac{5-4/(5\pi^2)}{3-4/(3\pi^2)}$, therefore, for $\forall i \in \mathbb{Z}^+$,
\begin{equation} \nonumber
\begin{aligned}
& \Big(4i-1-\dfrac{4}{(4i-1)\pi^2}\Big)\beta \geq \Big(4i+1-\dfrac{4}{(4i+1)\pi^2}\Big).
\end{aligned}
\end{equation}
Then by multiplying with $T_{bw}$, we have,
\begin{equation} \label{eq_t3_cond2}
\tilde{\delta}\tau_{i,1}^{max} \geq \tilde{\delta}\tau_{i,n_f}^{min}, \ \ \forall i \in \mathbb{Z}^+. 
\end{equation}
Lastly, we show that for any filter set $L_k(p), k = \overline{1,n_f}$ defined as in Theorem \ref{theorem_nf} with $n_f \geq (1/M+\beta-2)/(1/M-1), n_f \in \mathbb{N}$,
\vspace{-0.2cm}
\begin{equation} \nonumber
\forall \delta\tau_0 \neq 0, \exists L_q(p): \xi(L_q(p),\delta\tau_0) < \delta\tau_0.
\end{equation}
\vspace{0.3cm} \\
\textbf{\textit{Case 2-A:}} 
\begin{equation} \label{eq_t3_case2A}
\tilde{\delta}\tau_{i_0-1,1}^{min} < \delta\tau_0 \leq \tilde{\delta}\tau_{i_0,1}^{max}
\end{equation}
From (\ref{eq_t3_exlocprop}) and (\ref{eq_t3_case2A}), it is obvious that $\xi(L_1(p),\delta\tau_0) = \tilde{\delta}\tau_{i_0-1,1}^{min}$, and $\tilde{\delta}\tau_{i_0-1,1}^{min} < \delta\tau_0 $, hence $\xi(L_1(p),\delta\tau_0) < \delta\tau_0$.
\vspace{-0.3cm}
\begin{flushright}
$\square$
\end{flushright}
\vspace{-0.3cm}
\textbf{\textit{Case 2-B:}} 
\begin{equation} \label{eq_t3_caseb}
\tilde{\delta}\tau_{i_0-1,1}^{max} < \delta\tau_0 \leq \tilde{\delta}\tau_{i_0-1,1}^{min}.
\end{equation}
From (\ref{eq_t3_exloc}) we see that, $ \forall i \in \mathbb{Z}^+ $,
\begin{equation} \label{eq_t3_4}
\tilde{\delta}\tau_{i,n_f}^{min} < \tilde{\delta}\tau_{i,n_f-1}^{min} < ... < \tilde{\delta}\tau_{i,2}^{min} < \tilde{\delta}\tau_{i,1}^{min},
\end{equation}
hence,
\begin{equation} \label{eq_t3_5}
\tilde{\delta}\tau_{i_0-1,n_f}^{min} < \tilde{\delta}\tau_{i_0-1,n_f-1}^{min} < ... < \tilde{\delta}\tau_{i_0-1,1}^{min}.
\end{equation}
Now combining (\ref{eq_t3_cond2}) and (\ref{eq_t3_caseb}),
\begin{equation}  \nonumber
\begin{aligned}
& \tilde{\delta}\tau_{i_0-1,n_f}^{min} \leq \tilde{\delta}\tau_{i_0-1,1}^{max}, \\
&  \tilde{\delta}\tau_{i_0-1,1}^{max} < \delta\tau_0  \leq \tilde{\delta}\tau_{i_0-1,1}^{min},
\end{aligned}
\end{equation}
we have,
\begin{equation} \label{eq_t3_6}
\tilde{\delta}\tau_{i_0-1,n_f}^{min} < \delta\tau_0 < \tilde{\delta}\tau_{i_0-1,1}^{min}.
\end{equation}
Considering (\ref{eq_t3_5}) and (\ref{eq_t3_6}), there always exists a value $q\geq 2$ such that $\tilde{\delta}\tau_{i_0-1,q}^{min} < \delta\tau_0 < \tilde{\delta}\tau_{i_0-1,q-1}^{min}$. Following from (\ref{eq_t3_cond1}), we have $ \tilde{\delta}\tau_{i_0-1,q-1}^{min} \leq \tilde{\delta}\tau_{i_0,q}^{max}$. Therefore, 
\begin{equation} \label{eq_t3_7}
\tilde{\delta}\tau_{i_0-1,q}^{min} < \delta\tau_0  < \tilde{\delta}\tau_{i_0,q}^{max}.
\end{equation}
Combining (\ref{eq_t3_exlocprop}) and (\ref{eq_t3_7}),
\begin{equation}
\xi(L_{q}(p),\delta\tau_0) = \tilde{\delta}\tau_{i_0-1,q}^{min},
\end{equation}
therefore, $\xi(L_{q}(p),\delta\tau_0) < \delta\tau_0$. \hfill $\square$
\end{proof}

\end{document}